\documentclass[3p,11pt]{elsarticle}
\usepackage{amsmath}
\usepackage{graphicx}
\usepackage{caption} 
\usepackage{wrapfig}
\usepackage{color,soul}
\usepackage{enumitem}
\usepackage{caption}
\usepackage{subcaption}
\usepackage{comment}
\usepackage{amsmath}
\usepackage{times}
\usepackage{amssymb}
\usepackage{array}
\usepackage{latexsym}
\usepackage{amsthm}
\usepackage{amsfonts}
\usepackage{nomencl}
\usepackage[table]{xcolor}
\usepackage{lineno}

\makenomenclature
\usepackage{hyperref}
\usepackage{array,tabularx}
\usepackage{scalerel,stackengine}
\usepackage{graphics}
\usepackage[english]{babel}
\allowdisplaybreaks
\newtheorem{theorem}{Theorem}[section]

\newtheorem{lemma}[theorem]{Lemma}

\stackMath
\numberwithin{equation}{section}
\date{}
\newtheorem{df}{Definition}[section]

\usepackage[utf8]{inputenc}
\usepackage{algorithm}
\usepackage{algpseudocode}

\usepackage{multirow}
\usepackage{booktabs} 
\usepackage{array}    

\usepackage{amsmath, amssymb}
\usepackage{graphicx}
\usepackage{algorithm}
\usepackage{algpseudocode}
\usepackage{tikz}

\usetikzlibrary{shapes.geometric, arrows.meta, positioning}

\begin{document}
\begin{frontmatter} 
\title{Robust sensor coverage in the presence of spatial obstacles and exclusion zones\tnoteref{mytitlenote}} 
\author[1]{Vanshika Datta}
\author[1]{C. Nahak}
\address[1]{Department of Mathematics, Indian Institute of Technology Kharagpur, Kharagpur, West Bengal, India-721302} 

\cortext[mycorrespondingauthor]{Corresponding author} 
\ead{cnahak@maths.iitkgp.ac.in} 

\begin{abstract}
In this article, we develop a robust optimization framework for obstacle-aware aerial directional sensor networks operating under positional uncertainty using a sector-based sensing model. The monitoring area is discretized into grid points and a unified geometric formulation is developed to model directional sensing while explicitly accounting for obstacle-induced visibility loss and exclusion regions. To enhance sensing performance under uncertainty, three progressive optimization strategies, namely robust aerial grid coverage (RAGC), robust aerial target coverage (RATC) and robust aerial target scheduling (RATS), are proposed. The framework employs robust orientation optimization based on the radius of robust feasibility to maximize effective grid and target coverage under sensor perturbations, followed by a target-aware sleep scheduling mechanism that minimizes the number of active sensors without compromising target monitoring. Extensive simulations under varying target distributions, obstacle configurations, uncertainty levels and sensor failure scenarios demonstrate that the proposed framework achieves robust and energy-efficient sensing while consistently outperforming representative approaches from the literature in terms of coverage quality, target monitoring and sensor utilization.
\end{abstract}

\end{frontmatter}

\section{Introduction}
Wireless sensor networks have experienced substantial growth due to advances in micro-electromechanical systems (MEMS), wireless communication technologies and embedded electronics. Originally developed for military applications, WSNs are now widely used in environmental monitoring, healthcare, industrial automation, disaster management, precision agriculture and smart-city systems, see \cite{14,18,19}. In such applications, it is essential to ensure reliable monitoring of the region of interest and efficient service delivery, see \cite{7,9}. Among various quality-of-service (QoS) requirements, coverage is one of the most critical issues, as it directly influences both the sensing performance and deployment cost of the network. Consequently, sensing coverage has become a fundamental research problem in the design and development of WSNs in recent years.

Depending on the monitoring objective and the characteristics of the region of interest, coverage can be defined in several different ways. The three most widely studied coverage paradigms are area coverage, target coverage and barrier coverage. Area coverage aims to maximize sensing quality uniformly across the entire monitoring region and is commonly employed in applications such as forest surveillance, agricultural sensing and smart-city management \cite{3,12}. Barrier coverage seeks to establish a sensing boundary capable of detecting any object traversing a protected region and is particularly relevant to border surveillance, intrusion detection and military perimeter security \cite{11,19}. Target coverage focuses on monitoring a predefined set of critical points or objects while satisfying prescribed coverage requirements and is frequently applied to the protection of strategic facilities, transportation hubs, industrial plants and other critical infrastructure \cite{8,10}.

Beyond coverage objectives, WSNs can also be categorized according to the sensing capabilities of their deployed nodes. Homogeneous networks consist of sensors with identical sensing ranges, communication characteristics and energy resources, which simplifies network design and optimization. In contrast, heterogeneous networks incorporate sensors with diverse capabilities to enhance flexibility and monitoring effectiveness in complex environments \cite{13}. The assessment and optimization of coverage further depend on the spatial representation adopted for the monitoring region. Grid-based approaches discretize the area into finite cells and evaluate coverage at the cell level, whereas Voronoi-based methods partition the region into influence zones and are frequently used in distributed coverage optimization \cite{1,4,19}. In addition to achieving adequate coverage, energy efficiency remains a fundamental QoS requirement because sensor nodes typically operate under limited battery resources. Consequently, considerable research has focused on sleep scheduling, duty cycling, node selection and energy-aware deployment strategies to extend network lifetime while preserving the desired level of coverage \cite{2,14}.

Despite these advances, WSNs are frequently deployed in harsh and complex environments, requiring random sensor deployment, where coverage performance is affected by obstacle occlusion, blind spots, irregular boundaries and sensor malfunctions. Such factors alter the sensing geometry and can substantially reduce coverage quality and monitoring reliability. Moreover, practical deployments are inherently subject to uncertainty arising from localization errors, environmental disturbances, hardware imperfections, node failures and limited energy resources, making reliable coverage maintenance a challenging task \cite{24,4,15}. Furthermore, since sensor nodes are fragile devices susceptible to positional uncertainties, failures and energy constraints, robust coverage strategies are essential for maintaining reliable target monitoring in complex terrains. To address these issues, a variety of probabilistic, stochastic, fuzzy and robust optimization frameworks have been developed for coverage analysis and optimization \cite{12,21}. In particular, robust optimization has been widely adopted to account for uncertainty by seeking solutions that remain feasible under perturbations, thereby improving network reliability and resilience \cite{2,13}.

Robust optimization provides a systematic framework for incorporating uncertainty directly into the optimization process and obtaining solutions that remain feasible under a prescribed range of perturbations \cite{ro_book}. Owing to its ability to account for uncertainty without requiring detailed probabilistic information, robust optimization has been increasingly adopted in the design and analysis of engineering systems operating in uncertain environments \cite{ro_survey}. Nevertheless, assessing the degree of robustness associated with a feasible solution remains a challenging task. To address this issue, a quantitative measure, radius of robust feasibility (RRF), was introduced that characterizes the largest admissible uncertainty level for which feasibility is preserved \cite{rrf_start}. The exact computation and theoretical characterization of RRF have subsequently been established through rigorous optimization-based formulations \cite{rrf_exact}. More recently, the concept has been extended to broader optimization settings \cite{rrf_survey,rrf_end}. These developments have significantly enhanced its applicability as a practical tool for robustness assessment and uncertainty management in complex decision-making problems. For sensor networks, such a measure is particularly valuable because it provides an explicit robustness margin against deployment errors, environmental disturbances, obstacle interactions and other sources of uncertainty that can adversely affect coverage performance. These considerations become even more critical in directional sensor networks (DSNs), whose sensing behavior differs fundamentally from that of conventional omnidirectional sensor networks \cite{5,20}.

\begin{table}
\centering
\caption{Comparison of various QoS strategies}
\begin{tabular}{lllllll}
\hline
\textbf{QoS strategy} & \textbf{Sensing model} & \textbf{Occlusions} & \textbf{Lifetime} & \textbf{Sensor model} & \textbf{K-Target} & \textbf{Uncertainty} \\
\hline
PGAC \cite{1} & Sector & No & No & DSN & No & No \\ 
RO \cite{2} & Sector & No & Yes & WSN & No & Yes \\ 
EPSR \cite{5} & Grid & Yes & No & ASN & Yes & Partially \\
RESPIRE++ \cite{4} & Grid & No & No & DSN & Partially & Yes \\
IDS \cite{3} & Sector & No & No & DSN & No & No \\ 
Our & Grid & Yes & Yes & ASN & Yes & Yes \\
\hline
\end{tabular}
\label{tab:1}
\end{table}
Unlike conventional omnidirectional sensing models, DSNs employ sensors that monitor only a limited sector of the environment, making them suitable for applications involving cameras, radar systems and infrared sensors \cite{10,22}. The directional nature of sensing introduces additional decision variables and constraints, requiring the simultaneous consideration of sensing orientation, coverage overlap, blind regions and obstacle-induced visibility effects. As a result, coverage optimization in DSNs is substantially more challenging than in traditional omnidirectional networks. A distinctive feature of many directional sensors is their ability to adjust sensing orientation according to operational requirements, a capability commonly referred to as motility. This property has been extensively exploited to improve coverage quality and mitigate the adverse effects of overlap and occlusion. Nevertheless, existing studies have largely concentrated on coverage enhancement alone, while important QoS requirements such as target (k)-coverage, positional uncertainty, energy efficiency and network lifetime are often considered separately rather than within a unified framework \cite{16,17}.

Although significant progress has been made in DSN coverage optimization, existing methods often simplify sensing environments and rarely address obstacle-induced occlusions, irregular visibility constraints, positional uncertainty and energy-efficient operation within a unified framework. Consequently, robust coverage formulations that jointly consider obstacle awareness, exclusion zones, directional sensing, uncertainty handling and target coverage remain limited. Table~\ref{tab:1} summarizes representative studies, showing that most existing approaches satisfy only a subset of these QoS requirements. Furthermore, although obstacle-aware sensing has been investigated, the geometric visibility loss caused by obstacle occlusions is rarely integrated into the optimization process in a mathematically consistent manner. In particular, the effective sensing region obtained after obstacle clipping is seldom connected with computational coverage optimization, creating a disconnect between geometric modelling and practical deployment strategies. Accordingly, this study addresses several key questions: how obstacle-induced occlusions can be effectively incorporated into sensor coverage models; how environmental and deployment perturbations influence QoS performance and can be mitigated; and how multiple requirements, including deployment cost, positional uncertainty, target coverage and network lifetime, can be jointly optimized. Motivated by these challenges, the present study develops a unified robust coverage framework integrating obstacle-aware sensing, exclusion zones, uncertainty handling through RRF and energy-conscious optimization.

In this article, we develop a unified robust optimization framework for obstacle-aware target coverage in aerial directional sensor networks based on the RRF. Unlike existing approaches, which typically address obstacle avoidance, deployment uncertainty, target coverage, or energy efficiency independently, the proposed framework simultaneously integrates all these aspects into a single mathematical formulation. The sensing process is first modeled using a grid- and sector-based representation that explicitly captures obstacle-induced occlusions, exclusion zones, irregular sensing boundaries and directional sensing constraints, thereby providing a more realistic characterization of practical deployment environments. The exact obstacle-aware sensing geometry is then formulated analytically to determine the effective visible sensing region through geometric clipping of the nominal sensing sector. Since repeated analytical evaluation of these irregular clipped regions is computationally expensive during optimization, the proposed framework subsequently transforms the effective sensing region into a grid-based representation, enabling efficient numerical coverage optimization while preserving consistency with the underlying continuous geometric model. Sensor position uncertainty is subsequently incorporated through the RRF, enabling the sensing regions to remain reliable under deployment perturbations while excluding obstructed and non-coverable regions from the effective sensing footprint. Building upon this robust sensing model, three progressively enhanced optimization strategies are proposed to address different yet complementary design objectives. The robust aerial grid coverage (RAGC) strategy maximizes the effective sensing area, the robust aerial target coverage (RATC) strategy prioritizes reliable monitoring of target locations and the robust aerial target scheduling (RATS) strategy further introduces a target-aware sleep scheduling mechanism to significantly reduce the number of active sensors while maintaining high monitoring performance. Consequently, the proposed framework not only improves sensing robustness against environmental uncertainty and sensor failures but also achieves an effective balance between coverage quality, target monitoring and network lifetime. Extensive simulation studies demonstrate that the proposed methods consistently outperform representative approaches from the literature under varying deployment conditions, obstacle configurations, uncertainty levels and sensor failure scenarios. Overall, this work provides a mathematically rigorous and practically applicable framework for resilient, obstacle-aware and energy-efficient coverage optimization, thereby addressing several important limitations of existing robust coverage methodologies.

The remainder of the article is organized as follows. Section 2 introduces the necessary preliminaries and the basic sensing model. The proposed QoS-enhancing framework and the robust coverage model for spatial obstacles and exclusion zones are presented in Section 3. Section 4 details the proposed methodology and a solution algorithm. Simulation studies and performance comparisons with baseline methods are reported in Section 5. Finally, Section 6 concludes the paper and discusses possible directions for future research.

\section{Preliminaries}
This section introduces the mathematical concepts, sensing models and assumptions required for the proposed robust coverage framework. The principal notations and acronyms used throughout the paper are summarized in Table~\ref{tab:notation}. The required basic parameters, uncertainty representation, robust feasibility concepts, sensing geometries and associated mathematical tools are first presented, followed by the problem formulation and solution methodology.
\begin{table}
\centering
\caption{List of Notations}
\begin{tabular}{ll}
\hline
\textbf{Symbol} & \textbf{Description} \\
\hline
$m$ & Total number of sensors in the deployment. \\
$\mathcal{X}$ & Region of Interest.\\
$s$ & A particular sensor for reference. \\
$\mathcal{S} = \left\{ s^o_1, \cdots, s^o_m \right\}$ & Nominal location set of sensors. \\
$s^o_i=(x^o_i,y^o_i)$ & Nominal location of the $i^{th}$ sensor in $\mathbb{R} ^2$. \\
$s^w_i$ & Worst case location of  the $i^{th}$ sensor in $\mathbb{R} ^2$. \\
$s^{\rho}_i$ & Robust case location of the $i^{th}$ sensor in $\mathbb{R} ^2$. \\
$o_i$ & Orientation of the $i^{th}$ sensor after ground projection. \\
$\theta_s$ & Angle of view of the sensor. \\
$V_i$ & Selected vertex of orientation of the $i^{th}$ sensor. \\
$A_i(s_i,o_i)$ & Set of all covered points by the $i^{th}$ sensor in $\mathbb{R} ^2$. \\
$r_s$ & Minimum sensing range of a particular sensor $s$. \\
$\mathcal{U}_i^{\alpha}$ & Uncertainty set describing admissible perturbations in $\mathbb{R}^2$. \\
$\mathbb{U}_i^{\alpha}$ & Uncertainty set describing admissible perturbations in $\mathbb{R}^3$.\\
$\alpha$ & Radius of the uncertainty set. \\
$\rho_{i}$ & RRF associated with the $i^{th}$ sensor. \\
$VC(s_i)$ & Voronoi cell of the $i^{th}$ sensor.\\
$C_i{(o_i)}$ & Ground coverage of the $i^{th}$ sensor.\\
\hline
\end{tabular}
\label{tab:notation}
\end{table}

We use $\mathbb{R}^n$ to denote the $n$-dimensional Euclidean space, $||\cdot||$ for the Euclidean norm and $\mathbb{B}_n$ for the open unit ball in $\mathbb{R}^n$. For any vector $a$, its transpose is denoted by $a^T$. The Euclidean distance between two points $x,y \in \mathbb{R}^n$ is given by $d(x,y)=||x-y||$.

\subsection{Radius of robust feasibility} 
Consider the following parametric linear system under uncertainty in its constraint data, denoted by $\sigma^{\alpha}$:
\begin{equation}
    \sigma ^{\alpha } := \left\{ a^T_ix \leq b_i \right\}; \;(a_i,b_i) \in \mathcal{U} _i ^{\alpha} \subset \mathbb{R} ^{n+1}; \;i=\left\{1, 2, \cdots, p \right\} \;,
\end{equation}
where $(a_i,b_i) $ for $i=\left\{1, 2, \cdots, p \right\}$ are uncertain vectors and $\mathcal{U} _i ^{\alpha};\; i=\left\{1, 2, \cdots, p \right\}$ are the uncertain regions, typically assumed to be compact and convex.\\
The robust counterpart (RC) is obtained by ensuring that the constraints remain satisfied for all possible realizations within the predefined uncertainty set. Hence, the RC (deterministic form) of the corresponding system $\sigma ^{\alpha} $ with $U_i^{\alpha}=(\bar{ a_i}, \bar{b_i} ) + \alpha B_{n+1}, \; i=1:p$ and the corresponding feasible set $F_R^ {\alpha}$ for some $\alpha \in \mathbb{R} $ is:
\begin{equation}
    \sigma _ R ^ {\alpha} := \left\{ a_i^Tx \leq b_i ,\;  \forall (a_i,b_i) \in U_i^{\alpha},\; i=1:p \right\}.
\end{equation}
Robust optimization provides a framework for obtaining solutions that remain feasible under uncertain conditions. Instead of considering only nominal parameter values, robust formulations incorporate uncertainty directly into the decision-making process.

For a feasible solution, the radius of robust feasibility (RRF) quantifies the maximum uncertainty size for which the system constraints remain satisfied. The concept is defined as:
\begin{df}(Radius of robust feasibility)\label{rrf} 
Consider a parametric linear system subject to data uncertainty, denoted by $\sigma^{\alpha}$. Let $F_R^{\alpha}$ represent the feasible set of the robust counterpart of $\sigma^{\alpha}$. Then, the RRF for the system is given by: 
\begin{equation}
     \rho = \sup \{ \alpha \in \mathbb{R}_+ : (F_R ^ \alpha) \text{ is nonempty} \}.
\end{equation} 
\end{df}
The Minkowski functional is used to characterize the scaling of uncertainty sets. For a given convex set, the corresponding Minkowski functional and its properties are defined as:
\begin{df}(Minkowski function)
Let $\omega \subset \mathbb{R}^{n}$ be a convex set with $0_n$ in its interior. The Minkowski (or gauge) function associated with $\omega$, denoted by $\phi_{\omega}$, is the mapping $\phi_{\omega} : \mathbb{R}^{n} \rightarrow \mathbb{R}_{+}$, where $\mathbb{R}_{+} := [0,+\infty)$ and is defined as:
    \begin{equation}
        \phi_{\omega} (x) := \inf{ \left\{ t>0 : x \in t \omega \right\} },\; x \in \mathbb{R} ^n.
    \end{equation}
\end{df}
The following lemma summarizes the fundamental properties of the Minkowski function.
\begin{lemma}[See \cite{rrf_exact}, Lemma 2.2] Let $\omega \subset \mathbb{R}^{n}$ be a convex set whose interior contains $0_n$. Then, the following properties hold:
\begin{itemize}
    \item[$(1)$] $\phi_{ \omega}$ is sublinear and continuous.
    \item[$(2)$] $\left\{ x \in \mathbb{R} ^n : \phi _{\omega} (x) \leq 1 \right\} = cl (\omega)$, where $cl (\omega ) $ denotes the closure of $\omega$.
    \item[$(3)$] If $\omega$ is bounded and symmetric, $\phi _ {\omega }: = || \cdot ||$ defines a norm on $\mathbb{R} ^n$ generated by $\omega$.
\end{itemize} 
\end{lemma}
The exact expression for the RRF of an arbitrary set $S \subset \mathbb{R}^{n+1}$ was derived in the literature using the hypographical set $H(\bar{a},\bar{b})$ associated with the nominal system $\sigma$:
\begin{equation}
    H(\bar a, \bar b) := conv \left\{ (- \bar{a_i} , -\bar{b_i}): i=1:p \right\} + \mathbb{R} _+ (0_n, -1),
\end{equation}
where $\bar a = ( \bar{a_1},\bar{a_2}, \cdots, \bar{a_p} ) \in (\mathbb{R} ^n)^p \text{ and } \bar b = ( \bar{b_1},\bar{b_2}, \cdots, \bar{b_p} ) \in \mathbb{R} ^p$.
\begin{theorem} 
    If the nominal system is feasible, then the RRF of $\sigma^{\alpha}$ is given by:
    \begin{equation}
        \rho = \inf_{(a,b) \in H(\bar{a},\bar{b})} \phi_Z (a,b),
    \end{equation}
    where $Z\subset\mathbb{R}^{n}$ is a convex and compact set containing $0_n$ in its interior \cite{rrf_exact}. 
\end{theorem}
The obtained robustness radius provides a quantitative measure of tolerance against perturbations and is later utilized for sensor orientation adjustment.

\subsection{Sensing model in 2D}
Coverage performance depends strongly on the sensing characteristics of deployed sensors. Two sensing models are generally considered in the literature: omnidirectional sensor networks (OSNs) and directional sensor networks (DSNs). In OSNs, sensors provide uniform sensing capability in all directions within a predefined sensing radius. Therefore, the sensing region is generally represented as a circular area. In contrast, DSNs employ directional sensing mechanisms where the sensing region is restricted by both sensing range and angular field of view. Such sensors include camera-based, radar and infrared sensing systems.

\begin{figure}
    \centering
    \includegraphics[width=0.5\textwidth]{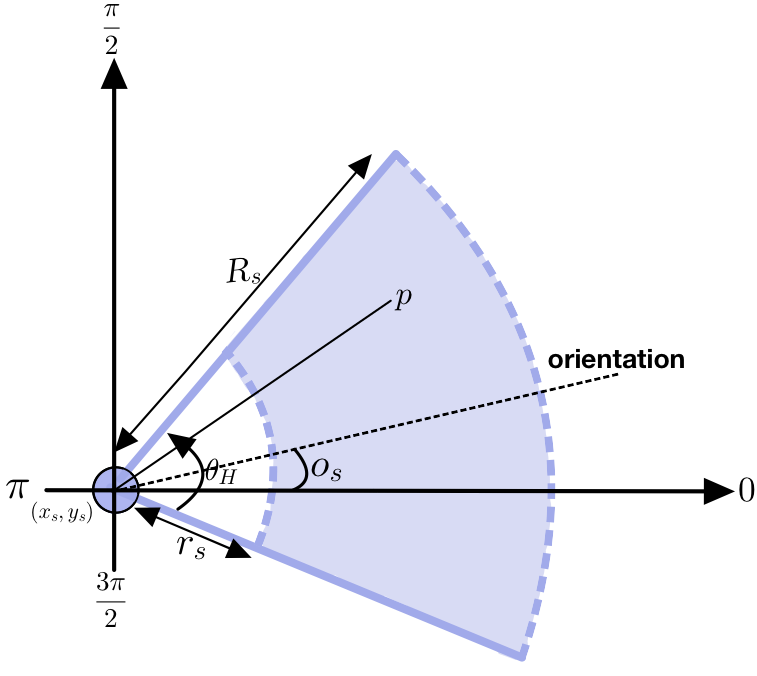}
    \caption{Two-dimensional projection of the aerial directional sensing model.}
    \label{Fig_3}
\end{figure}
Let $\mathcal{X} \subset \mathbb{R}^{2}$ denote the monitoring region and $\mathcal{S}$ represent the set of deployed sensors with a particular sensor $s \in \mathcal{S}$. In this work, each directional sensor is represented on the ground plane by its location $s=(x_s,y_s)\in\mathbb{R}^{2}$. The sensing capability of a directional sensor is characterized by its sensing range and angular field of view, see Figure~\ref{Fig_3}. Let $r_s$ and $R_s$ denote the minimum and maximum sensing radii, respectively and let $\theta_H$ represent the horizontal field of view. The orientation of sensor $s$ is defined by the angle $o_s$, measured with respect to the positive $x$-axis. The sensing region of a directional sensor is represented as an annular sector on the ground plane. For a point $p=(x,y)\in\mathbb{R}^{2}$, the point is considered covered by sensor $s$ if it satisfies both the radial and angular constraints. The effective sensing region $C_s(o_s)$ is defined as:
\begin{equation}
\label{sensing_area}
C_s(o_s)=
\left\{
(x,y)\in\mathbb{R}^{2}\,\;\middle|\,\;
r_s\leq ||(x_s,y_s)-(x,y)||\leq R_s,\,\,
\operatorname{atan2}(y-y_s,x-x_s)
\in
\left[
o_s-\frac{\theta_H}{2},\,
o_s+\frac{\theta_H}{2}
\right]
\right\}.
\end{equation}
Here, the angle between the vector from the sensor $s$ to the point $(x,y)$ and the positive $x$-axis is computed using $\operatorname{atan2}(y-y_s,x-x_s)$. Therefore, the coverage region of each directional sensor is determined by the parameters $(r_s,R_s,\theta_H,o_s)$.

\subsection{Problem overview and assumptions}
This work considers the robust target $k$-coverage problem for an aerial directional sensor network deployed over a region of interest $\mathcal{X},$ containing obstacles and exclusion zones. Each sensor provides a sector-shaped sensing region, whose effective coverage may be affected by positional uncertainty, orientation errors, visibility constraints and possible sensor failures. The objective is to determine the sensor orientations that maximize the effectively covered ground area while preserving coverage feasibility under bounded uncertainty.

To reduce redundant coverage evaluation and assign local sensing responsibility, $\mathcal{X}$ is discretized into a finite set of grid points. Then we represent a few possible target and obstacles location and the coverage status is determined based on the visibility and sensing constraints of the deployed sensors. This discretization converts the continuous coverage problem into a finite optimization problem while allowing obstacle-induced coverage loss and irregular sensing boundaries to be incorporated explicitly.

The following assumptions are considered:
\begin{itemize}
    \item Sensor locations are known with bounded positional uncertainty.
    \item Obstacles and exclusion zones modify the sensing regions by blocking or removing coverage areas.
    \item Sensor orientations are adjustable within their allowable angular ranges.
    \item Sensor failures and redundant sensing are considered during coverage optimization.
    \item The uncertainty level is represented using a predefined uncertainty set for robust feasibility analysis.
\end{itemize}

The proposed methodology first determines the obstacle-aware feasible sensing regions by removing obstructed and non-coverable areas. Subsequently, the RRF-based robust orientation strategy is employed to maintain coverage under perturbations, followed by an energy-aware mechanism to deactivate redundant sensors while satisfying the required coverage level.

\section{Model formulation}
In this section, we introduce the sensing system used throughout the paper, along with its functionalities and control parameters. Then, we define two critical surveillance requirements, i.e., field of view and image resolution. Finally, we convert the real-life problem with obstacles into a mathematical optimization problem and formulate the proposed model.

\subsection{Sensing model in 3D}
\begin{figure}
    \centering
    \includegraphics[width=0.8\textwidth]{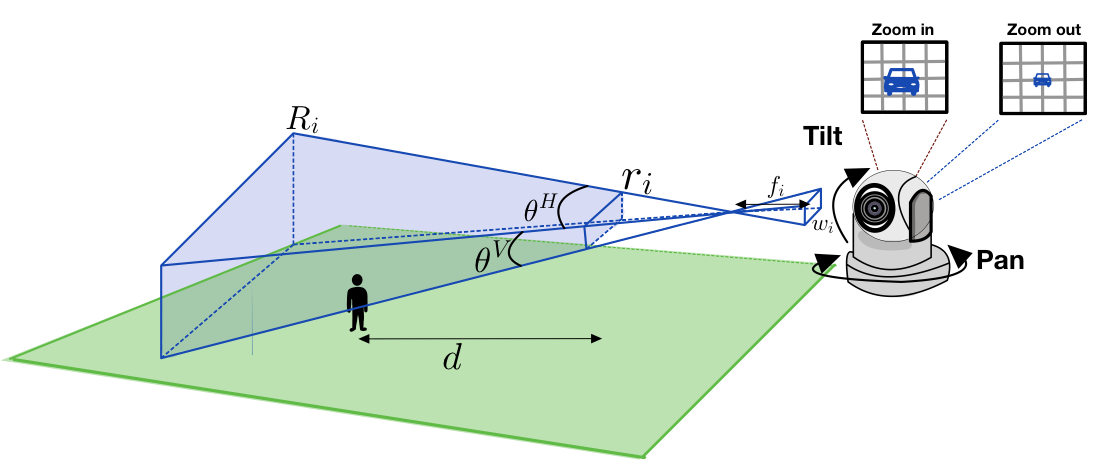}
    \caption{PTZ camera architecture and the corresponding three-dimensional directional sensing model.}
    \label{Fig_4}
\end{figure}

Unlike conventional omnidirectional sensors, DSNs possess the ability to dynamically adjust their sensing direction and field of view according to monitoring requirements. For example, pan-tilt-zoom (PTZ) cameras possess horizontal, vertical and zoom capabilities, as shown in Figure \ref{Fig_4}. Such sensors are widely employed in surveillance, environmental monitoring and target-tracking applications due to their flexibility in reducing blind regions and improving target visibility.

Consider a set of $m$ sensors with location set $$S=\left\{ s_1, s_2, \cdots, s_m \right\}$$ with the $i^{th}$ sensor located at $s_i=(x_i,y_i,z_i).$ Each sensor is characterized by three controllable parameters: the pan angle $\theta ^H _i$, the tilt angle $\theta ^V _i$ and the zoom parameter $f _i$, denoting its focal length. The pan angle governs horizontal rotation about the vertical axis, whereas the tilt angle controls vertical orientation. The zoom parameter determines the camera magnification and consequently influences the effective sensing region. Hence, the sensing region of each sensor can be characterized by its field of view $\mathcal{F } _i = \left\{ ( \theta ^V  _i, \theta ^H  _i, f _i ) \right\}.$
The pan angle is related to the focal length and width of the image sensor, say $w_i$ and the tilt angle is related to the focal length and length of the image sensor, say $l_i$. This leads to the relations
\begin{equation}
        \operatorname{tan}{\left(\frac{\theta ^H  _i }{2}\right)}  = \left( \frac{w_i}{2 f_i} \right)
        \;\; \text{and}\; \;
        \operatorname{tan}{\left(\frac{\theta ^V  _i }{2}\right)} = \left( \frac{l_i}{2 f_i} \right).
\end{equation}
Next, the values of $r_i$ and $R_i$ are to be determined. Observe that $r_i$ is equal to half of the pan focus distance, i.e., $2r_i=f_i$ \cite{r_s}. Another important parameter is the imaging resolution, which signifies the number of pixels corresponding to each unit foot of an object in the image. It is one of the crucial factor in surveillance QoS and is given by:
$$\text{Resolution}= \frac{n_h}{\mathcal{F}_h},$$
where, $n_h$ represents the number of horizontal pixels of the sensor image and $\mathcal{F } _h$ represents the horizontal field, i.e. the width of the object at a distance $d$, see Fig~\ref{Fig_4}. It follows that
$$\mathcal{F } _h = d \cdot \frac{w_i}{f_i}.$$
Referring to Figure~\ref{Fig_4}, the two triangles corresponding to the horizontal field of view of the sensor are similar. Therefore, the ratio between $R_i$ and $\mathcal{F}_h$ is equal to that between $f_i$ and $w_i$. Multiplying both sides by $n_h$ yields
\begin{equation*}
    \begin{split}
    &\frac{R_i}{\mathcal{F } _h} \cdot n_h = \frac{f_i}{w_i} \cdot n_h \\
    \implies & R_i = \frac{f_i}{w_i} \cdot \frac{n_h}{n_h} \cdot \mathcal{F } _h\\
    \implies & R_i = \frac{f_i}{w_i} \cdot \frac{n_h}{\text{Resolution}} \cdot
    \end{split} 
\end{equation*}
Here, $R_i$ denotes the distance between the pan focal distance and the image sensor positioned at infinity, commonly referred to as the hyperfocal distance. Accordingly, surveillance is feasible only when the distance between the sensor and the target lies within the interval $(r_i,,R_i)$, where $r_i$ and $R_i$ represent the minimum and maximum permissible distances, respectively \cite{R_s}.

\begin{figure}
    \centering
    \includegraphics[width=\textwidth]{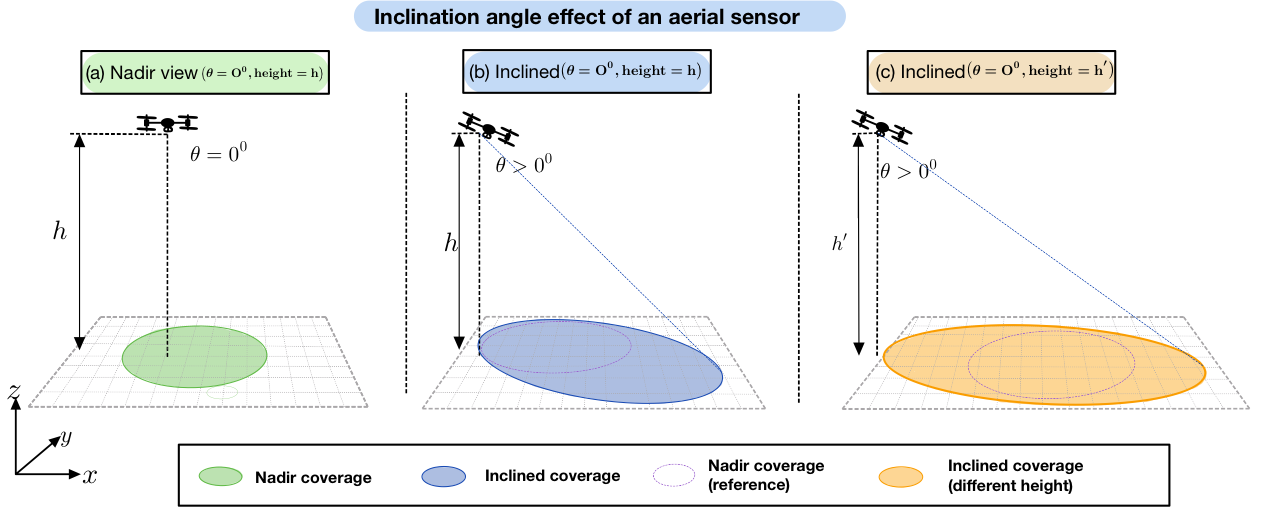}
    \caption{Effect of sensor inclination and altitude on the projected ground coverage.}
    \label{Fig_5}
\end{figure}

The coverage region of a directional sensor depends on its position, sensing radius, orientation angle and field-of-view parameter. Consequently, variations in sensor inclination directly influence the projected coverage area. As the inclination angle changes, the sensing footprint undergoes geometric deformation, altering effective coverage and visibility, see Figure~\ref{Fig_5}. Note that when the view angle, i.e., the orientation of the sensor, is larger, it indicates that the object's face is more inclined in the captured image. The more the inclination, the less clear is the image recognition \cite{img_recognition}. The zoom functionality modifies the field of view and sensing characteristics of the camera. In general, larger zoom values provide greater target resolution while reducing the observable area. 

One of the key advantages of PTZ sensors is their orientation adaptability. By adjusting the pan and tilt angles, the sensor can redirect its sensing sector toward regions of interest and compensate for visibility degradation caused by obstacles or environmental changes. Let $\Delta\theta ^H _i$ and $\Delta\theta ^V _i$ denote the pan and tilt adjustments, respectively. The updated orientation is therefore given by
\[
\theta_{i}^{H\text{new}}=\theta_{i}^{H}+\Delta\theta_{i}^{H},
\qquad \text{and} \qquad
\theta_{i}^{V\text{new}}=\theta_{i}^{V}+\Delta\theta_{i}^{V}.
\]
These adjustable orientation parameters provide the foundation for the robust coverage optimization framework proposed in this study, where sensor directions are continuously refined to improve target coverage in the presence of uncertainty, obstacles and exclusion zones.

\subsection{Surveillance constraint}
We can verify whether a sensor has a direct line of sight to a target object by checking for an intersection between the orientation vector of the sensor and the surface of the physical obstacles. Here, several existing technologies can convert these physical obstacles in images into coordinate descriptions and their surface representations, see \cite{los1,los2,los3}.

Consider the $i^{th}$ sensor located at $(x_i,y_i,z_i)$, whose sensing direction is represented by the unit orientation vector $\mathbf{o}_i$. For an arbitrary target point $t_j=(x_j,y_j,z_j)$, the relative position vector from the sensor to the target is defined as
$$\mathbf{r}_{ij} = (x_j-x_i,\;y_j-y_i,\;z_j-z_i).$$ 
Moreover, the Euclidean distance between the sensor and target is computed as 
$$d(s_i,t_j)=\sqrt{(x_j-x_i)^2+(y_j-y_i)^2+(z_j-z_i)^2}.$$ 
A target is considered detectable only if it lies within the sensing range and the field of view of the PTZ sensor. Let $R_i$ and $r_i$ denote the maximum and minimum sensing range, respectively and $\theta ^V _i$ denote the field-of-view angle. Then the distance constraint is
$$r_i \leq d(s_i,t_j) \leq R_i.$$
This relation indicates that the object should be able to meet the requirements of imaging resolution. Hence, if it satisfies this equation, it indicates that the target is within the depth distance of $\mathcal{F} _h$. Let us denote this coverage condition by $\tau_{ij}$. Then it takes the value 1 if the target is within reach, else 0 as follows:
\begin{equation} \label{constraint1}
\tau_{ij}= 
\begin{cases}
    1, & \text{if } r_i \leq d(s_i,t_j) \leq R_i\\
    0, & \text{otherwise. } 
\end{cases}
\end{equation}
Furthermore, the angular deviations in horizontal and vertical directions, say $\delta ^H $ and $\delta ^V$, between the sensor orientation and the target direction must satisfy the following:
\begin{equation} \label{Eqn1}
    \begin{split}
        o_i - \frac{\theta ^H _i}{2} \leq \delta ^H = \cos^{-1}\left(\frac{\mathbf{o}_i \cdot \mathbf{r}_{ij}}{\|\mathbf{r}_{ij}\|}\right) ^ H \leq o_i + \frac{\theta ^H _i}{2},\\
        o_i - \frac{\theta ^V _i}{2} \leq \delta ^V = \cos^{-1}\left(\frac{\mathbf{o}_i \cdot \mathbf{r}_{ij}}{\|\mathbf{r}_{ij}\|}\right) ^ V \leq o_i + \frac{\theta ^V _i}{2}.
    \end{split}
\end{equation}
The above condition ensures that the target lies within the sensing sector view angle generated by the PTZ sensor. The coverage condition, denoted by $\alpha_{ij}$, is defined as follows:
\begin{equation} \label{constraint2}
\alpha_{ij}= 
\begin{cases}
    1, & \text{if satisfying Eqn.}~\eqref{Eqn1} \\
    0, & \text{otherwise. } 
\end{cases}
\end{equation}
Finally, introduce $\beta^{req}$ as the minimum required angle to ensure that the inclination angle $\theta_{ij}$ between the $i^{th}$ sensor and the $j^{th}$ target does not become too large for effective surveillance. Accordingly, we consider the following:
\begin{equation} \label{constraint3}
\beta_{ij}= 
\begin{cases}
    1, & \text{if } \theta_{ij} \leq \beta ^{ req}  \\
    0, & \text{otherwise. } 
\end{cases}
\end{equation}
The surveillance constraint is expressed as
\begin{equation}
    \gamma _{ij} = \tau _{ij} \times \alpha _{ij} \times \beta _{ij}, \;\; \gamma _{ij} \in \left\{ 0,1 \right\},
\end{equation}
where $\gamma _{ij} $ indicates whether the $j^{th} $ target can be covered by the $i ^{th}$ sensor or not.

\subsection{Modelling obstacles}
\begin{figure}
\centering
\begin{subfigure}{0.36\textwidth}
    \centering
    \includegraphics[width=\linewidth]{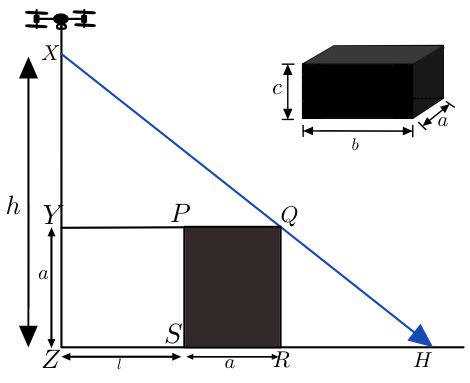}
    \caption{Side view of the occlusion region generated by the obstacle's height.}
\end{subfigure}
\hfill
\begin{subfigure}{0.36\textwidth}
    \centering
    \includegraphics[width=\linewidth]{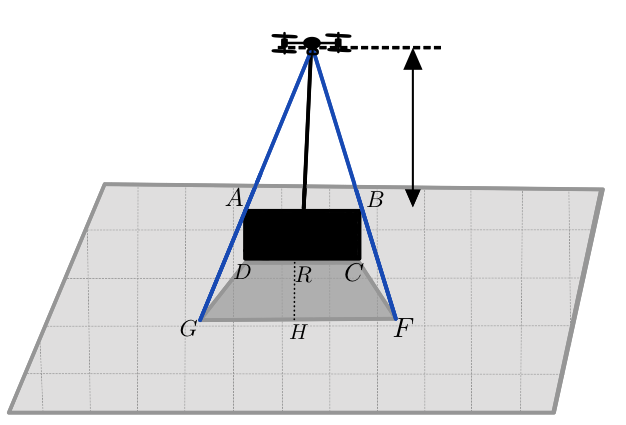}
    \caption{Top view of the occlusion region generated by the obstacle length.}
\end{subfigure}
\vspace{0.4cm}
\begin{subfigure}{0.50\textwidth}
    \centering
    \includegraphics[width=\linewidth]{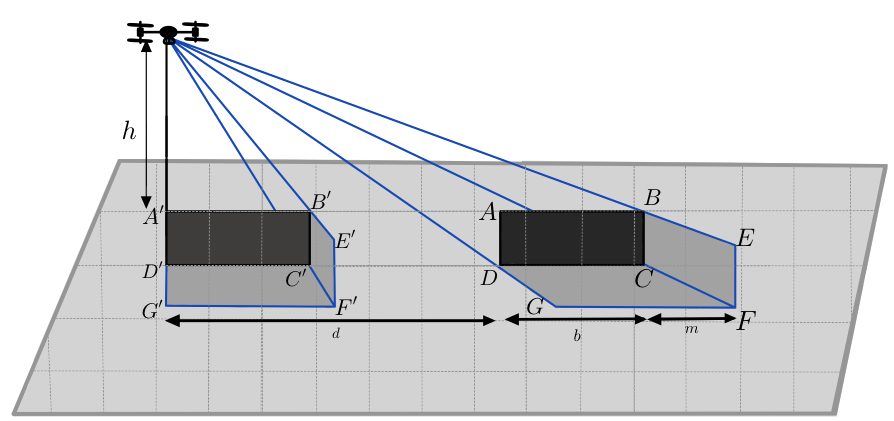}
    \caption{Geometric representation of the occlusion region caused by the obstacle.}
\end{subfigure}
\caption{Illustration of obstacle-induced occlusion regions in side and top views together with the corresponding geometric model.}
\label{Fig_6}
\end{figure}

Consider an obstacle with dimensions $a \times b \times c$ in three-dimensional space. The corresponding surveillance region occluded by the obstacle can be divided into two parts, as illustrated in Figure~\ref{Fig_6}. The first occluded region, denoted by $O_1$, arises due to the shelter of the obstacle  DCFG, as shown in Figure~\ref{Fig_6b}. The second region, denoted by $O_2$, is caused by the longer part BCFE, as depicted in Figure~\ref{Fig_6c}. Hence, the total occluded area is obtained by determining these two components separately.

Based on the geometric relationship, the first part of the occluded area $O_1$ can be deduced as: 
\begin{align}
        O_1 & = O_{DCFG}= \frac{(|DC| + |GF|)}{2} \cdot |RH|, \label{O1_1}\\
        \frac{|AB|}{|GF|} &=\frac{|XY|}{|XZ|} = \frac{(h-c)}{h}, \label{O1_2} \\
        \frac{|ZH|}{|RH|} &= \frac{|XZ|}{|QR|} = \frac{h}{c}. \label{O1_3}
\end{align}
Using \eqref{O1_2} and the fact that $AB=b$, corresponding to the longer dimension of the obstacle, the equation can be rewritten as:
\begin{equation}
    |GF|= \left( \frac{bh}{h-c} \right). \label{gf}
\end{equation}
Let $l$ denote the distance between the obstacle and the pole on which the sensor is deployed; see Figure~\ref{Fig_6a}. Substituting this distance into \eqref{O1_2} yields
\begin{equation}
\begin{split}
\frac{a+l+|RH|}{|RH|} = \frac{h}{c}
\implies |RH|= \frac{(a+l) \cdot c}{(h-c)}. \label{rh}
\end{split}
\end{equation}
Combining \eqref{O1_1}, \eqref{gf} and \eqref{rh}, we obtain
\begin{equation*}
O_1 = \left( |DC| + \frac{bh}{h-c} \right) \left( \frac{(a+l) \cdot c}{2 \cdot (h-c)} \right) = \frac{(bc) \cdot (a+l) \cdot (2h-c)}{2 \cdot (h-c)^2}. \label{O_1}
\end{equation*}
Similarly, for the second part of the occluded region, we have the following relations: 
\begin{equation}\label{O2_}
    \begin{split}
        O_2 & = O_{BCFE}= \frac{(|BC| + |FE|)}{2} \cdot m, \\
        \frac{|CB|}{|FE|} &=\frac{|B'C'|}{|F'E'|} = \frac{|XY|}{|XZ|}= \frac{(h-c)}{h},  \\
        \frac{|D'C|}{|G'F|} & = \frac{b+d}{b+d+m}. 
    \end{split}
\end{equation}
Solving \eqref{O2_} gives:
$$O_2 = \frac{(ac) \cdot (2h-c) }{2 \cdot (h-c)^2} \cdot (b+d).$$
Combining this expression with \eqref{O_1}, the total occluded area is given by
\begin{equation*}
   \begin{split}
        O_{area} &= O_1 + O_2\\
    &=  \frac{c(2h-c)(ad+bl+2ab)}{2 (h-c)^2}.
   \end{split}
\end{equation*}
For obstacles with irregular geometries where closed-form occlusion calculation becomes difficult, the obstacle boundary can be represented using finite line segments and the feasible sensing region can be obtained through geometric intersection checks. The effective sensing footprint of each sensor is accordingly clipped by obstacle-induced visibility constraints and combined with location uncertainty through a robustified framework, as we shall see in the next subsection.

\begin{figure}
\begin{subfigure}{\textwidth}
  \centering
  \includegraphics[width=\linewidth]{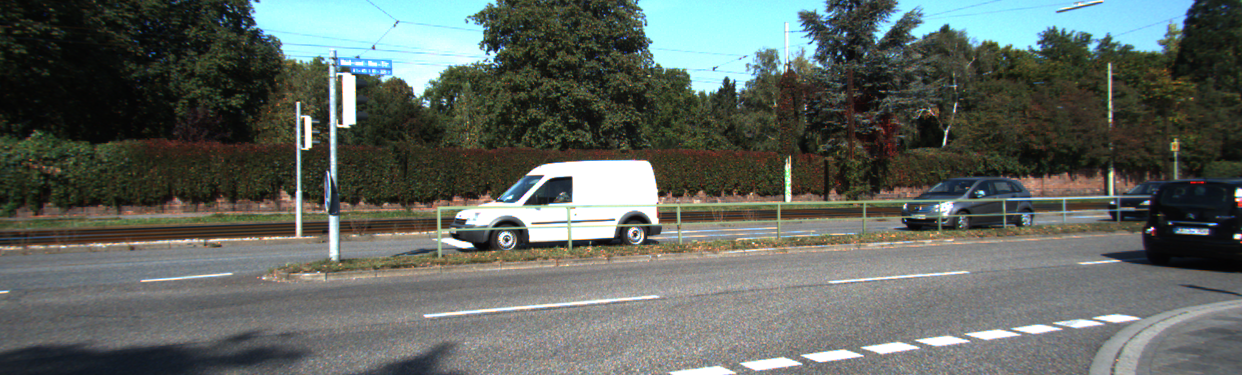}
  \caption{Real-world sensing scenario illustration.}
  \label{Fig_3a}
\end{subfigure}%

\begin{subfigure}{.5\textwidth}
  \centering
  \includegraphics[width=1.2\linewidth]{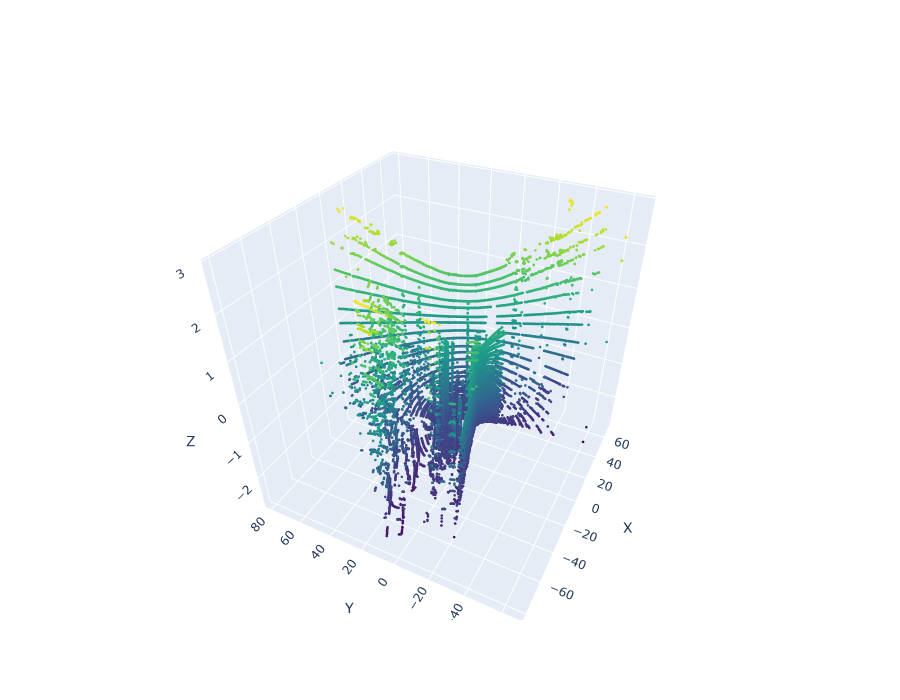}
  \caption{LiDar 3D point cloud.}
  \label{Fig_3b}
\end{subfigure}
\begin{subfigure}{.5\textwidth}
  \centering
  \includegraphics[width=\linewidth]{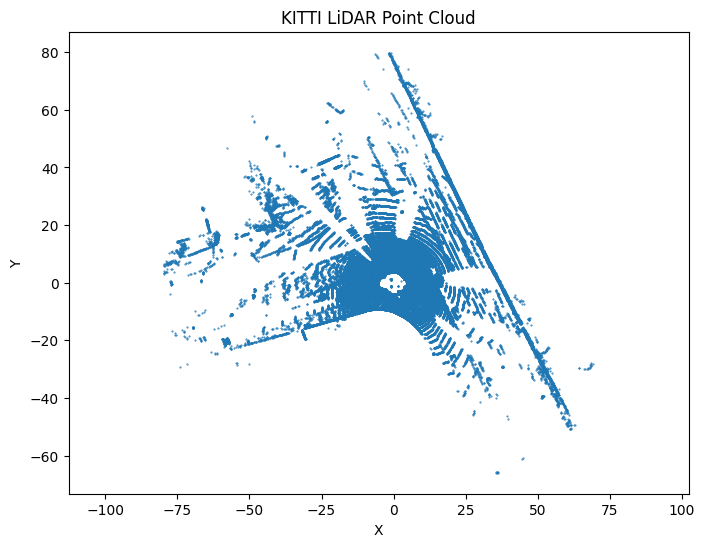}
  \caption{Bird eye view.}
  \label{Fig_3b}
\end{subfigure}
\caption{Example urban traffic scene from the KITTI dataset: (a) RGB image showing vehicles, road structures and occlusions; (b) corresponding LiDAR point cloud providing 3D spatial information; (c)  bird-eye view representation in 2D. Images are obtained from the KITTI Vision Benchmark Suite.}
\label{Fig_3ab}
\end{figure}

Figure~\ref{Fig_3ab} illustrates a representative RGB image of an urban driving scenario from the KITTI dataset and the corresponding LiDAR point cloud representation of the same traffic scene, see \cite{KITTI}. Although simulation platforms such as CARLA Simulator provide flexible and controllable environments for generating synthetic sensor data, the KITTI dataset is employed in this work to illustrate real-world driving scenarios captured under actual operating conditions and the potential line-of-sight obstructions. The LiDAR representation illustrates the role of accurate 3D sensing in environments affected by visibility constraints, while the proposed framework remains applicable to both real-world and simulated scenarios. 

\subsection{Modelling uncertainty and problem formulation}
In practical deployments, the nominal sensor locations may not represent the actual deployed positions, even when sensors are placed according to a predefined configuration or randomly distributed. Such positional deviations can significantly influence the sensing regions and may lead to inaccurate orientation decisions, reduced coverage quality and violation of coverage requirements. Therefore, incorporating uncertainty into the optimization process is essential to ensure reliable network performance. In this work, the radius of robust feasibility (RRF) is employed to quantify the maximum uncertainty level within which the coverage constraints remain feasible, enabling robust sensor orientation optimization and improved quality of service. 

We consider a set of $m$ sensors with nominal location set $$S = \left\{ s^o_1,s^o_2, \cdots, s^o_m \right\}\subset \mathbb{R} ^2,$$ where $s^o_i=(x^o_i,y^o_i)$ denotes the ground-projected nominal location of the $i^{th}$ sensor for $i=1:m$. To capture deployment inaccuracies, each sensor location is allowed to deviate from its nominal position ${s^o_i}$ within the uncertainty set:
\begin{align}\label{uncertainity}
    \mathcal{U}_i^{\alpha} = {s^o_i} + \alpha \mathbb{B}_2,
\end{align}
where the scalar $\alpha>0$ represents the radius of the uncertainty set and quantifies the maximum admissible deviation of the actual sensor location from its nominal position.

Let $\mathcal{X}$ represent the region of interest, which is discretized into a finite set of grid points. The objective is to maximize the overall coverage, measured by the total number of grid points monitored by the sensors. For each sensor location $s_i^o \in S$, the corresponding set of covered grid points is determined. Let $g$ denote the total number of grid points and define the variable $G_j \in \{0,1\}$ be such that
\begin{equation} \label{grid_point}
    G_j = 
    \begin{cases}
        1, & \text{if } j^{th} \text{ grid point is covered} \\
        0, & \text{otherwise. } 
    \end{cases}
\end{equation}
For each element in $S$, we compute the number of grid points covered over all admissible locations within its uncertainty set $\mathcal{U}_i^{\alpha}$. Accordingly, we consider the worst-case approximate location in sense that it covers minimum number of grid points in any orientations, $s^w_{i}$ of the $i^{th}$ sensor corresponding to its nominal position $s^o_i$ along its sensing direction $\vec{u}$ as
\begin{equation}
    s^w_{i} = {s^o_i} + \alpha \frac{\vec u}{|\vec u|}.
\end{equation}
Under this model, coverage must remain feasible for all realizations $s_i \in \mathcal{U}_i^{\alpha}$. To quantify this uncertainty tolerance, we employ the concept of RRF, originally introduced in \cite{12}. In this work, the RRF of a sensor denotes the largest perturbation radius for which the coverage constraints remain satisfied while maintaining a prescribed coverage threshold $\delta$. Let $\rho_{s_i}$ denote the RRF corresponding to the $i^{th}$ sensor and $\rho_{\min}$ denote a user-defined robustness requirement representing the minimum perturbation level that must be tolerated by the system. Enforcing this condition ensures that the obtained sensing configuration remains feasible under location deviations of at least $\rho_{\min}$.  Based on the computed RRF, we consider the following robust approximate location of the $i^{th}$ sensor along its sensing direction $\vec{u}$:
\begin{equation}
    s^{\rho}_{i} = {s^o_i} + \rho_i \frac{\vec u}{|\vec u|},
\end{equation}
which provides a tractable representation for evaluating coverage under admissible perturbations. 

Building on this idea, the resulting optimization problem is formulated as
\begin{equation} 
\text{(Coverage maximization):  } \max_{o_i} \sum_{j=1}^{g} G_j(s^{\rho}_{i},o_i),\;i=1:m;\; j=1:g;
\end{equation}
where, $G_j$ is given by equation~\eqref{grid_point}. Since this problem is a mixed-integer linear programming problem (MILP), i.e., $X_i \in \left\{ 0,1 \right\};\; \theta^H_i, \theta^V_i \in [-\pi,\pi];\; f_i\in \mathbb{R}, $ it is classified as NP-hard \cite{NP_hard}. This formulation provides the orientations of sensors that maximize the number of grid points covered, hence the coverage. 

\section{Proposed methodology}
In this section, we present the proposed methodology and the algorithm developed to solve the formulated robust optimization problem. The effective visible sensing region is first obtained by geometrically clipping the nominal sensing sector according to obstacle-induced occlusions and exclusion zones. This clipped sensing region is subsequently discretized into uniformly spaced grid points, providing a computationally efficient representation for coverage evaluation and optimization while preserving the underlying geometric visibility model. Sensor position uncertainty is incorporated through the RRF to construct robust feasible sensing regions. Building upon this unified geometric-computational framework, three progressively enhanced optimization strategies, namely robust aerial grid coverage (RAGC), robust aerial target coverage (RATC) and robust aerial target scheduling (RATS), are introduced to maximize grid coverage, improve target coverage and minimize the number of active sensors through target-aware sleep scheduling, respectively. The complete solution procedure is summarized as a set of algorithms.

\subsection{Obstacle-aware sector decomposition and coverage evaluation}

\begin{figure}
\centering
\begin{subfigure}{0.45\columnwidth}
    \centering
    \includegraphics[width=\linewidth]{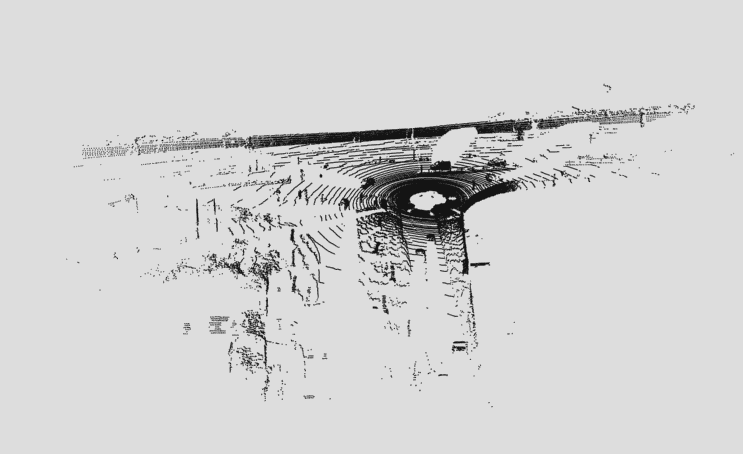}
    \caption{Point cloud representation of the environment.}
    \label{Fig_3d}
\end{subfigure}
\hfill
\begin{subfigure}{0.45\columnwidth}
    \centering
    \includegraphics[width=\linewidth]{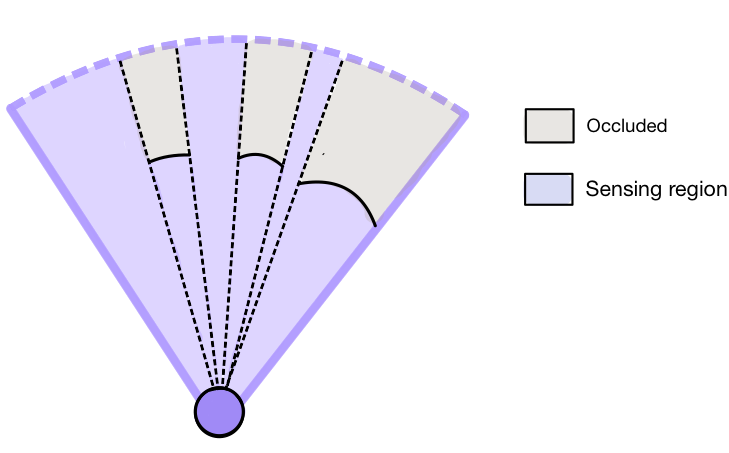}
    \caption{Effective sensing sector split after obstacle-induced occlusion.}
    \label{Fig_3e}
\end{subfigure}
\caption{Illustration of obstacle-induced occlusion in directional sensor coverage generated from the real scene shown in Figure~\ref{Fig_3}.}
\end{figure}

We first consider the point cloud data shown in Figure~\ref{Fig_3} and reconstruct it from a lower viewing angle. The right half of Figure~\ref{Fig_3d} corresponds to the reconstructed scene shown in Figure~\ref{Fig_3}, while Figure~\ref{Fig_3e} illustrates its sector-wise partition after removing the occluded regions.

In a DSN, the sensing region is represented as a truncated sector-based coverage area, as depicted in Figure~\ref{Fig_3e}. Due to the presence of obstacles and exclusion zones, the original sensing region is partitioned into multiple feasible sector regions. The objective is to maximize the effective sensing area, i.e., to maximize the number of covered grid points within the discretized region of interest $\mathcal{X}$. Therefore, it is necessary to identify the grid points covered by each feasible sensor sector and evaluate their contribution to the overall coverage.

To determine the coverage contribution of each feasible sector region, the intersection points between the grid structure and the sector boundaries are first identified. These boundary intersections define the limits of the covered region and enable the accurate determination of grid points within the sensing area. The coverage is then evaluated by checking whether the grid points lie within the resulting sector boundaries. For the intersection analysis, six possible configurations can occur, as illustrated in Figure~\ref{Fig_9}. For simplicity, a local coordinate system is adopted, with the sensor positioned at $(x,0)$ after suitable translation and rotation. Let $h$ represent the perpendicular distance of an intersection point from the $x$-axis. Each configuration represents a distinct geometric interaction between the grid points and the sector boundary, which is analyzed individually in the following cases:
\begin{figure}
\centering
\captionsetup{justification=centering}

\begin{subfigure}{0.32\textwidth}
  \centering
  \includegraphics[width=0.5\linewidth]{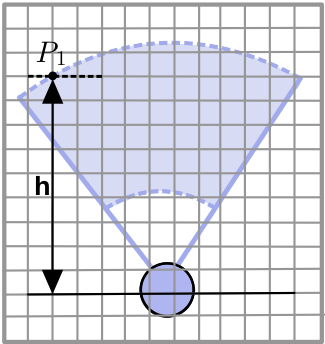}
  \caption{Case 1}
  \label{Fig_case1}
\end{subfigure}
\begin{subfigure}{0.32\textwidth}
  \centering
  \includegraphics[width=0.5\linewidth]{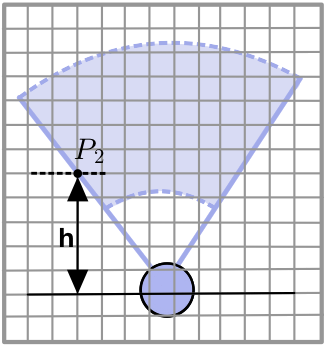}
  \caption{Case 2}
  \label{Fig_case2}
\end{subfigure}
\begin{subfigure}{0.32\textwidth}
  \centering
  \includegraphics[width=0.5\linewidth]{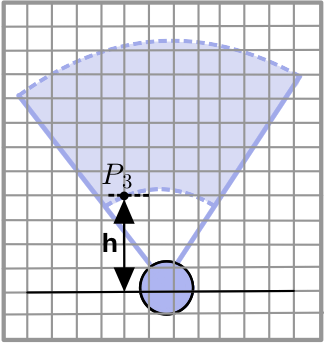}
  \caption{Case 3}
  \label{Fig_case3}
\end{subfigure}

\begin{subfigure}{0.32\textwidth}
  \centering
  \includegraphics[width=0.5\linewidth]{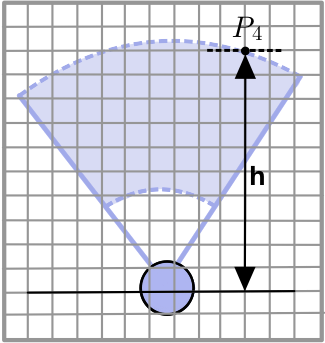}
  \caption{Case 4}
  \label{Fig_case4}
\end{subfigure}
\begin{subfigure}{0.32\textwidth}
  \centering
  \includegraphics[width=0.5\linewidth]{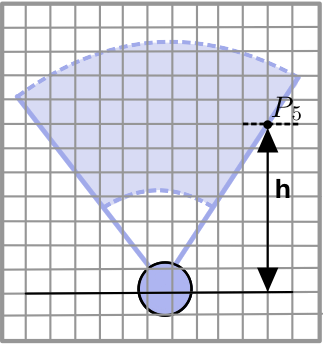}
  \caption{Case 5}
  \label{Fig_case5}
\end{subfigure}
\begin{subfigure}{0.32\textwidth}
  \centering
  \includegraphics[width=0.5\linewidth]{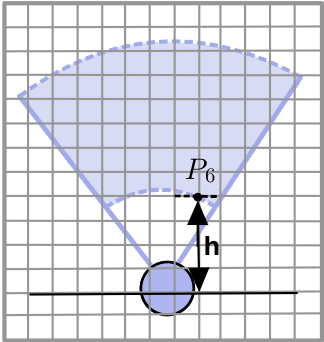}
  \caption{Case 6}
  \label{Fig_case6}
\end{subfigure}

\caption{Representative sensing sector and grid lines intersection configurations.}
\label{Fig_9}
\end{figure}

\subsubsection{Case 1}
This case corresponds to the configuration in which the grid intersects the left side of the outer arc of the sector boundary, as illustrated in Figure~\ref{Fig_case1}. When
$$\pi - 2 \cdot\arcsin{\frac{h}{R_s}} \leq \theta_s$$
and
$$\left(\arcsin{\frac{h}{R_s}}\right),\;
\left(\pi-\arcsin{\frac{h}{R_s}}\right)
\in
\left[
\left(\beta_s-\frac{\theta_s}{2}\right),
\left(\beta_s+\frac{\theta_s}{2}\right)
\right],$$
the intersection points lie only on the outer arc of the sensing region. The corresponding coordinate of $P_1$ is
$$
P_1=x_s-R_s\cos\left(\arcsin{\frac{h}{R_s}}\right).
$$

\subsubsection{Case 2}

In this configuration, the intersection between the Voronoi edge and the sensing region satisfies the following conditions:
$$\arcsin{\frac{h}{R_s}}, \pi - \arcsin{\frac{h}{R_s}} \notin \left[ \beta _s - \frac{\theta _s}{2}, \beta _s + \frac{\theta _s}{2}  \right],$$
$$\pi - 2\cdot \arcsin{\frac{h}{R_s}} >\theta _s \;\; \text{and}\; \;\arcsin{\frac{h}{R_s}} < \beta_s < \pi - \arcsin{\frac{h}{R_s}}.$$
Under these conditions, the Voronoi edge intersects the sector sidelines, as shown in Figure~\ref{Fig_case2}. Based on the location of the intersection points, the coordinate of $P_2$ is:
$$
P_2=\left(x_s  - \frac{h}{\tan{ \left( \beta_s + \frac{\theta_s}{2} \right) }} \right).
$$

\subsubsection{Case 3}

This configuration is characterized by the conditions
$$ \pi - \arcsin{\frac{h}{r_s}} \in \left[ \left( \beta_s - \frac{\theta_s}{2} \right), \left( \beta_s + \frac{\theta_s}{2} \right)  \right],$$ 
$$ \pi - \arcsin{\frac{h}{R_s}} \notin \left[ \left( \beta_s - \frac{\theta_s}{2} \right), \left( \beta_s + \frac{\theta_s}{2} \right)  \right],$$
and
$$ \pi - 2 \cdot \arcsin{\frac{h}{r_s}} > \theta_s,$$
as displayed in Figure~\ref{Fig_case3}. The corresponding coordinate of $P_3$ is
$$P_3= \left(x_s  - \frac{h}{\tan{ \left( \beta_s + \frac{\theta_s}{2} \right) }} \right).$$

\subsubsection{Case 4}

This case corresponds to the situation in which
$$\pi - 2 \cdot\arcsin{\frac{h}{R_s}} \leq \theta_s$$ and $$ \left(\arcsin{\frac{h}{R_s}}\right), \left(  \pi - \arcsin{\frac{h}{R_s}} \right) \in \left[ \left( \beta_s - \frac{\theta_s}{2} \right), \left( \beta_s + \frac{\theta_s}{2} \right)  \right].$$ 
Under these conditions, the intersection points are located exclusively on the outer arc of the sensing region (see Figure~\ref{Fig_case4}). The coordinate of $P_4$ is
$$ P_4 = x_s + R_s \cos{\left(\arcsin{\frac{h}{R_s}}\right)}.$$

\subsubsection{Case 5}

This configuration arises when the intersection between the Voronoi edge and the sensing region satisfies
$$\arcsin{\frac{h}{R_s}}, \pi - \arcsin{\frac{h}{R_s}} \notin \left[ \beta _s - \frac{\theta _s}{2}, \beta _s + \frac{\theta _s}{2}  \right],$$
$$\pi - 2\cdot \arcsin{\frac{h}{R_s}} >\theta _s \;\; \text{and}\; \;\arcsin{\frac{h}{R_s}} < \beta_s < \pi - \arcsin{\frac{h}{R_s}},$$ as demonstrate in Figure~\ref{Fig_case5}. In this scenario, the Voronoi edge intersects the sector sidelines. Accordingly, the coordinate of $P_5$ is
$$
P_5=\left(x_s  + \frac{h}{\tan{ \left( \beta_s - \frac{\theta_s}{2} \right) }} \right).
$$

\subsubsection{Case 6}

This case is identified by the conditions
$$ \arcsin{\frac{h}{r_s}} \in \left[ \left( \beta_s - \frac{\theta_s}{2} \right), \left( \beta_s + \frac{\theta_s}{2} \right)  \right],$$ 
$$ \arcsin{\frac{h}{R_s}} \notin \left[ \left( \beta_s - \frac{\theta_s}{2} \right), \left( \beta_s + \frac{\theta_s}{2} \right)  \right],$$ 
and 
$$ \pi - 2 \cdot \arcsin{\frac{h}{r_s}} > \theta_s,$$ 
as depicted in Figure~\ref{Fig_case6}. The coordinate of $P_6$ is
$$P_6=\left(x_s  + \frac{h}{\tan{ \left( \beta_s - \frac{\theta_s}{2} \right) }} \right).$$
\subsection{Analytical computation of effective coverage area}

After identifying the grid points lying inside the effective sensing region, the corresponding continuous coverage area can be obtained analytically using elementary geometric formulations. The analytical formulation presented below provides the exact geometric description of the effective sensing region and enables quantitative evaluation of the actual visible sensing area.

Let the effective sensing region be represented by the clipped sensing sector after removing the obstacle-induced occluded region. The boundary of this region is completely determined by the sector sidelines, circular arc boundaries, and the obstacle-intersection points obtained in the previous subsection. Consequently, the effective sensing region can be decomposed into a finite number of elementary geometric components, each corresponding to either a circular sector, a triangle, or a truncated sector.

The area of each triangular component is computed using Heron's formula

\begin{equation}
\mathrm{Area}
=
\sqrt{d(d-e_1)(d-e_2)(d-e_3)},
\end{equation}

where

\[
d=\frac{e_1+e_2+e_3}{2}
\]

denotes the semiperimeter and $e_1,e_2,e_3$ denote the side lengths of the corresponding triangle.

To determine the boundary of the visible sensing sector, the two sector sidelines are first obtained as

\begin{equation}
\begin{cases}
(a_0-x_s)\cos\beta_s
+
(b_0-y_s)\sin\beta_s
=
\sqrt{(a_0-x_s)^2+(b_0-y_s)^2}
\cos\left(\frac{\theta_H}{2}\right),\\
r_s
\le
d(s,p)
\le
R_s,
\end{cases}
\end{equation}

where $s=(x_s,y_s)$ denotes the sensor location, $\beta_s$ represents the sensor orientation, $\theta_H$ is the horizontal field of view, and $p=(a_0,b_0)$ denotes any point lying on the sideline.

Similarly, the inner and outer sector boundaries are determined using the corresponding circular arc equations

\begin{equation}
\begin{cases}
d(s,p)=r_s,\\
a_0
\in
\left[
r_s\cos\left(\beta_s+\frac{\theta_H}{2}\right)+x_s,
~
r_s\cos\left(\beta_s-\frac{\theta_H}{2}\right)+x_s
\right],\\
b_0
\in
\left[
r_s\sin\left(\beta_s-\frac{\theta_H}{2}\right)+y_s,
~
r_s\sin\left(\beta_s+\frac{\theta_H}{2}\right)+y_s
\right],
\end{cases}
\end{equation}

and

\begin{equation}
\begin{cases}
d(s,p)=R_s,\\
a_0
\in
\left[
R_s\cos\left(\beta_s+\frac{\theta_H}{2}\right)+x_s,
~
R_s\cos\left(\beta_s-\frac{\theta_H}{2}\right)+x_s
\right],\\
b_0
\in
\left[
R_s\sin\left(\beta_s-\frac{\theta_H}{2}\right)+y_s,
~
R_s\sin\left(\beta_s+\frac{\theta_H}{2}\right)+y_s
\right].
\end{cases}
\end{equation}

Using these boundary equations together with the obstacle-intersection points computed previously, the effective sensing region can be partitioned into elementary geometric components whose areas are evaluated individually and combined to obtain the exact visible sensing area corresponding to the sensor orientation.

Although the analytical formulation provides the exact effective sensing area, directly optimizing this continuous geometric region is computationally expensive because the visible sensing boundary changes irregularly with sensor orientation and obstacle configuration. Therefore, the region of interest is discretized into uniformly spaced grid points, and the effective sensing region is represented numerically by the subset of grid points lying inside the clipped sensing polygon. Consequently, maximizing the number of covered grid points serves as an efficient approximation of maximizing the effective visible sensing area while preserving the underlying analytical model.

\subsection{Robust aerial grid coverage optimization}

Based on the discretized representation of the effective sensing region, the objective of the first optimization strategy is to maximize the total number of grid points covered by the network while satisfying the sensing, obstacle, and robustness constraints. Let the discretized region of interest be represented by the finite set of grid points,
\[
\mathcal{G}=\{t_1,t_2,\ldots,t_n\},
\]
where each \(t_i\) denotes a grid point belonging to the feasible sensing region. Accordingly, the robust aerial grid coverage problem is formulated as
\begin{equation}
\max_{\{\beta_s\}_{s\in S}}
\sum_{p\in\mathcal{G}} I(p),
\end{equation}
where
\[
I(p)=
\begin{cases}
1,& \text{if } p \text{ belongs to the effective sensing region},\\
0,& \text{otherwise}.
\end{cases}
\]

The optimization is performed subject to the sensing constraints, obstacle-induced visibility constraints, exclusion zones, and the robust feasibility constraints. Since the resulting optimization problem is nonlinear and combinatorial, an iterative orientation optimization algorithm is developed to obtain a computationally efficient solution. The proposed RAGC algorithm is summarized in Algorithm~\ref{alg:RAGC}.

\begin{algorithm}
\caption{RAGC: Robust aerial grid coverage}
\label{alg:RAGC}
\begin{algorithmic}[1]

\Require Sensor set $\mathcal{S}$, grid points $\mathcal{G}$, obstacle set $\mathcal{O}$, RRF values
\Ensure Robust sensor orientations maximizing effective grid coverage

\For{each active sensor $s_i$}
    \State Initialize the best orientation $\beta_i^\ast$ using the current orientation.
    \For{each candidate orientation $\beta=0,\Delta\beta,\ldots,360-\Delta\beta$}
        \State Determine grid points satisfying sensing range, field-of-view and line-of-sight constraints.
        \State Remove obstacle and exclusion grid points.
        \State Compute the total effective coverage objective.
        \If{current objective exceeds the best objective}
            \State Update $\beta_i^\ast$.
        \EndIf
    \EndFor
    \State Assign $\beta_i\leftarrow\beta_i^\ast$.
\EndFor
\State Return optimized sensor orientations.

\end{algorithmic}
\end{algorithm}

\subsection{Robust target coverage and \texorpdfstring{$k$}{k}-coverage analysis
}

After obtaining the feasible sensing regions, the coverage status of each grid point in the discretized region of interest ($\mathcal{X}$) is evaluated. Let $\mathcal{G}=\{t_1,t_2,\dots,t_g\}$ denote the set of grid points and $G_{ij}$ indicate whether sensor $i$ covers grid point $t_j$, defined as
\begin{equation}
G_{ij}=
\begin{cases}
1, & t_j \text{ is covered by $i^{th}$ sensor},\\
0, & \text{otherwise}.
\end{cases}
\end{equation}
The coverage level of a grid point is calculated as
\begin{equation}
K_j=\sum_{i=1}^{m} G_{ij},\; j = 1: g,
\end{equation}
where $K_j$ represents the number of sensors covering the grid point $t_j$. For reliable monitoring applications, single sensor coverage may not be sufficient due to possible failures or uncertainties. Therefore, the framework is extended to the target $k$-coverage problem, where each target point is required to be covered by at least $k$ sensors. For target coverage, each required target point should satisfy $K_j\geq1$, whereas for target $k$-coverage, the constraint becomes
\begin{equation}
K_j \geq k, \quad \forall t_j\in \mathcal{G}.
\end{equation}
Therefore, the objective is to maximize the number of grid points satisfying the required coverage condition while considering obstacle-induced visibility restrictions and uncertainty.

To enhance the target monitoring capability of the proposed framework, the robust aerial target coverage (RATC) algorithm is introduced. It determines sensor orientations that maximize the number of covered targets while reducing unnecessary overlap among neighboring sensors whenever possible. The complete procedure is summarized in Algorithm~\ref{alg:RATC}.
\begin{algorithm}
\caption{RATC: Robust aerial target coverage}
\label{alg:RATC}
\begin{algorithmic}[1]

\Require Optimized RAGC solution, target set $\mathcal{T}$, obstacle set $\mathcal{O}$
\Ensure Robust sensor orientations maximizing target coverage

\For{each active sensor $s_i$}
    \State Initialize the best orientation $\beta_i^\ast$.
    \For{each candidate orientation $\beta$}
        \State Identify all visible targets satisfying sensing range, field-of-view and line-of-sight constraints.
        \State Remove targets already assigned to neighbouring sensors whenever possible.
        \State Compute the number of newly covered targets.
        \If{target coverage improves}
            \State Update $\beta_i^\ast$.
        \EndIf
    \EndFor
    \State Assign $\beta_i\leftarrow\beta_i^\ast$.
\EndFor
\State Return target-oriented sensor configuration.

\end{algorithmic}
\end{algorithm}

\subsection{Energy-efficient sensor scheduling}
After satisfying the required coverage condition, redundant sensors are identified to reduce unnecessary energy consumption. A binary activation variable $z_i$ is introduced for each sensor, where
\begin{equation}
z_i=
\begin{cases}
1, & \text{sensor } i \text{ remains active},\\
0, & \text{sensor } i \text{ enters sleep mode}.
\end{cases}
\end{equation}
The sleep scheduling process aims to minimize the number of active sensors while preserving the target $k$-coverage requirement. Hence, the following condition must remain satisfied:
\begin{equation}
\sum_{i=1}^{m} z_i G_{ij} \geq k, \quad \forall t_j\in\mathcal{G}.
\end{equation}
Thus, redundant sensors can be temporarily deactivated without affecting the required sensing reliability, improving energy efficiency and network lifetime.

To improve the energy efficiency of the proposed framework, the robust aerial target scheduling (RATS) algorithm is proposed. It identifies redundant sensors whose monitored targets remain covered by other active sensors and schedules them into sleep mode, while sensors covering no target are deactivated. The proposed procedure is outlined in Algorithm~\ref{alg:RATS}.

\begin{algorithm}
\caption{RATS: Robust aerial target scheduling}
\label{alg:RATS}
\begin{algorithmic}[1]

\Require Target-oriented sensor configuration obtained from RATC
\Ensure Energy-efficient active sensor configuration

\For{each sensor $s_i$}
    \State Determine the targets monitored by $s_i$.
    \If{$s_i$ monitors no target}
        \State Assign inactive state.
    \Else
        \State Check whether every monitored target is also covered by another active sensor.
        \If{all monitored targets remain protected}
            \State Schedule $s_i$ into sleep mode.
        \Else
            \State Keep $s_i$ active.
        \EndIf
    \EndIf
\EndFor
\State Return the final sensor orientations and active sensor set.

\end{algorithmic}
\end{algorithm}

\section{Simulation results and performance analysis}
This section presents the simulation-based evaluation of the proposed obstacle-aware robust coverage framework under uncertain sensor deployment. The simulation environment, deployment strategy and parameter settings are first described, followed by the visualization of the proposed orientation optimization and target-aware sleep scheduling procedures. Subsequently, the performance of the proposed framework is evaluated through quantitative measures, including grid coverage, target coverage and the number of active sensors. Finally, a comparative analysis with representative approaches from the literature is conducted to demonstrate the effectiveness, robustness and energy efficiency of the proposed framework under different deployment scenarios, obstacle configurations and sensor failure conditions.

\subsection{Simulation setup}
The proposed framework is implemented in Python 3.12 using the Spyder integrated development environment. All simulations are performed over a two-dimensional square region of interest, where aerial directional sensors are randomly deployed and the sensing field is discretized into uniformly spaced grid points. Rectangular obstacles are randomly generated to model occlusions and inaccessible regions, while target points are randomly selected from feasible grid locations outside the obstacle regions. To account for deployment uncertainty, the RRF is computed for each sensor, within which its position is randomly perturbed before applying the proposed orientation optimization strategies. Subsequently, a target-aware sleep scheduling mechanism is employed to reduce the number of active sensors while preserving effective target monitoring. The framework is evaluated under different deployment scenarios by varying the number of sensors, target points, obstacle configurations, grid resolution, sensing parameters and sensor failures. Performance is assessed in terms of effective grid coverage, target coverage and the number of active sensors and is compared with representative methods from the literature to demonstrate the performance of the proposed approach.

\subsection{Orientation optimization and target-aware sleep scheduling} 
Figure~\ref{fig:deployment_optimization} illustrates the progressive orientation optimization achieved by the proposed framework. The first figure shows the initial random deployment, while the remaining figures present the optimized sensor orientations obtained using RAGC, RATC and RATS, respectively. It can be observed that the optimized deployments significantly reduce unnecessary overlaps and improve sensing effectiveness by orienting sensors towards more informative regions despite obstacle-induced occlusions.

The quantitative comparison of the three proposed strategies is presented in Table~\ref{tab:our_comparison}. RAGC primarily maximizes the effective sensing area and therefore achieves the highest grid coverage of 699 grid points using 43 active sensors, although only 20 targets are monitored. By explicitly prioritizing target coverage, RATC nearly doubles the number of covered targets to 39, with an expected reduction in grid coverage to 611 grid points while retaining the same number of active sensors. Finally, RATS introduces target-aware sleep scheduling, reducing the number of active sensors from 43 to only 26 while maintaining coverage of 38 targets, resulting in only a marginal loss of one target despite a substantial reduction in active sensors. These results demonstrate that the three proposed strategies effectively address different deployment objectives, namely maximizing grid coverage, maximizing target coverage and achieving an energy-efficient target monitoring solution with minimal performance degradation.
\begin{table}
\centering
\caption{Comparison of the proposed optimization strategies}
\begin{tabular}{llll}
\hline
\textbf{Method} & \textbf{Active sensors} & \textbf{Grid points covered} & \textbf{Targets covered} \\
\hline
Initial & 50 & 526 & 17\\
RAGC & 43 & 699 & 20 \\
RATC & 43 & 611 & 39 \\
RATS & 26 & 383 & 38 \\
\hline
\end{tabular}
\label{tab:our_comparison}
\end{table}

\begin{figure}
\begin{subfigure}{0.5\textwidth}
  \centering
  \includegraphics[width=\linewidth]{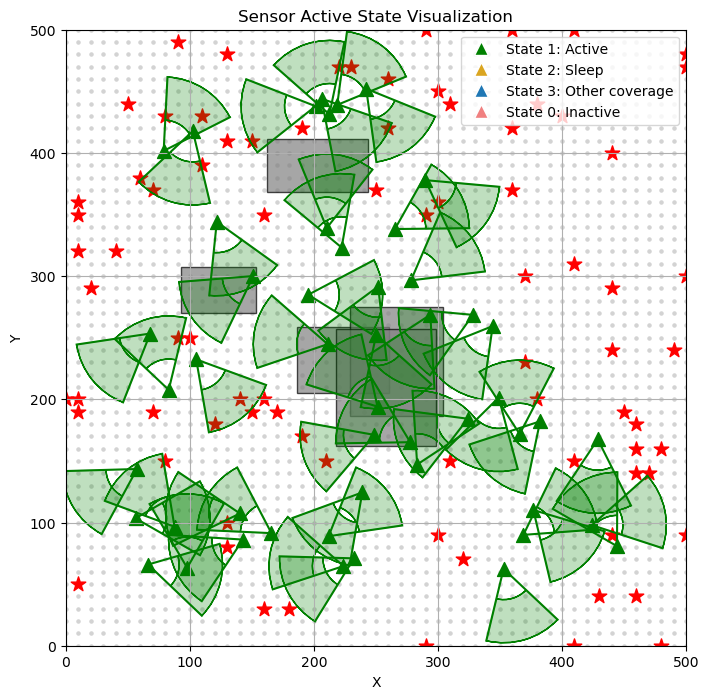}
  \caption{Initial sensor deployment.}
  \label{Fig_3b}
\end{subfigure}
\begin{subfigure}{0.5\textwidth}
  \centering
  \includegraphics[width=\linewidth]{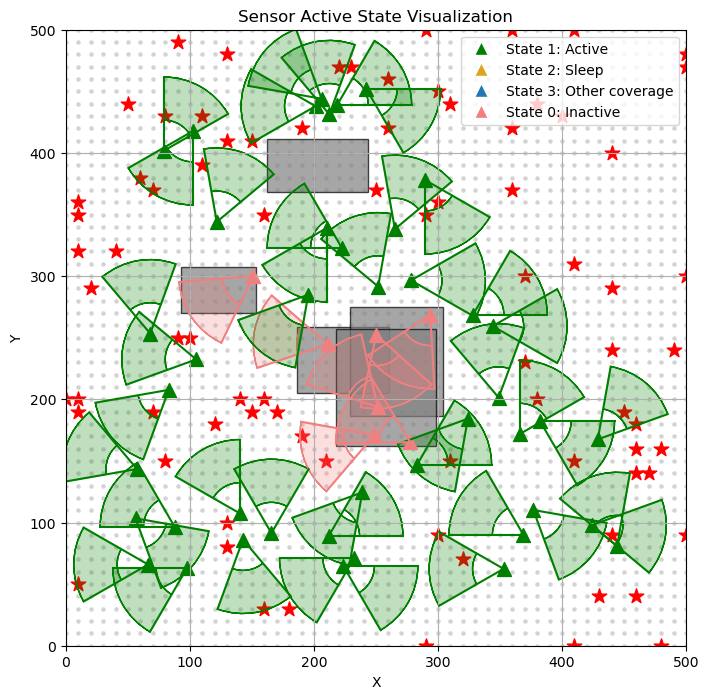}
  \caption{RAGC: Robust Aerial Grid Coverage.}
  \label{Fig_3b}
\end{subfigure}\\
\begin{subfigure}{0.5\textwidth}
  \centering
  \includegraphics[width=\linewidth]{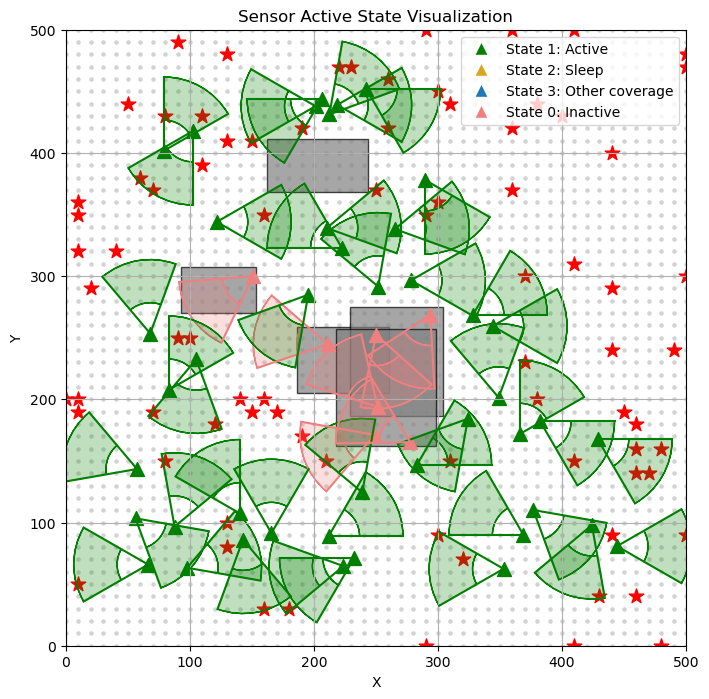}
  \caption{RATC: Robust Aerial Target Coverage.}
  \label{Fig_3b}
\end{subfigure}
\begin{subfigure}{0.5\textwidth}
  \centering
  \includegraphics[width=\linewidth]{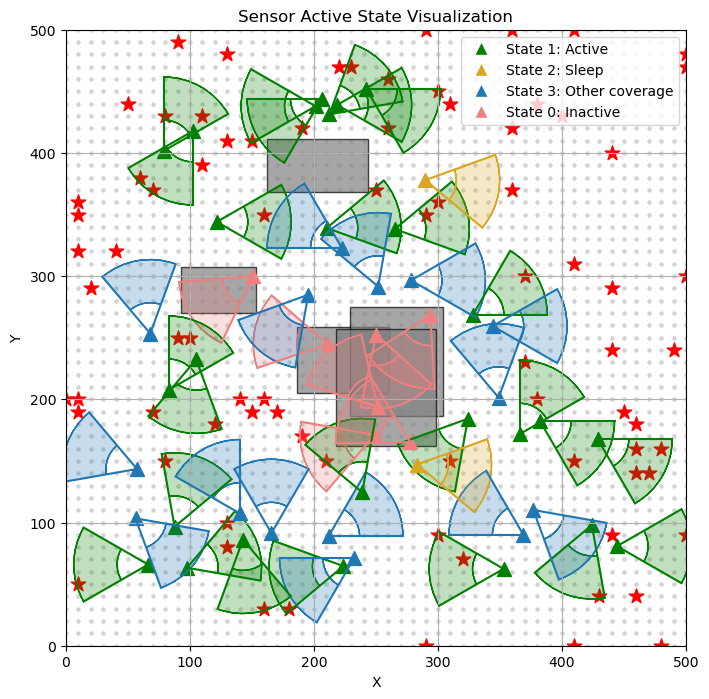}
  \caption{RATS: Robust Aerial Target Scheduling.}
  \label{Fig_3b}
\end{subfigure}
\caption{Progressive orientation optimization of the proposed framework.}
\label{fig:deployment_optimization}
\end{figure}

\subsection{Performance evaluation and comparative analysis}
To further evaluate the proposed framework, extensive simulations are conducted under different deployment settings by varying important network parameters. Since the proposed optimization framework is derived from the analytical effective sensing region, we first compare the effective sensing area obtained by different approaches before evaluating the corresponding grid and target coverage.
\begin{figure}
    \centering
    \includegraphics[width=0.6\textwidth]{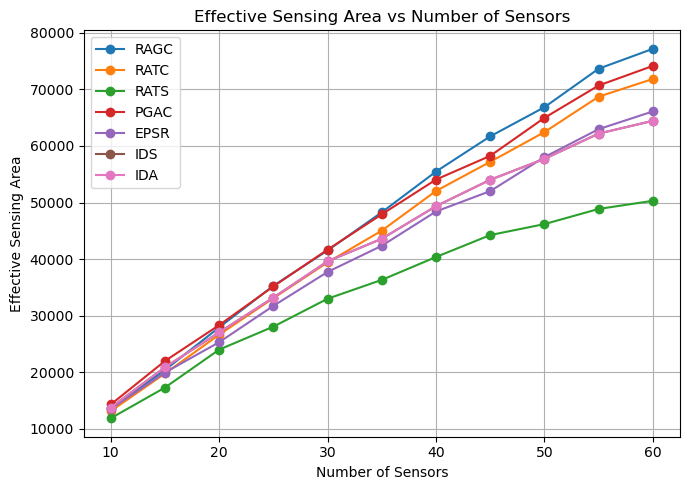}
    \caption{Variation of effective sensing area with respect to the number of deployed sensors for different optimization methods.}
    \label{fig:effective_area}
\end{figure}
\begin{table}
\centering
\caption{Effective sensing area obtained by different optimization methods under varying numbers of deployed sensors.}
\label{tab:effective_area}
\resizebox{\textwidth}{!}{
\begin{tabular}{lccccccccccc}
\hline
\textbf{Method} & \textbf{10} & \textbf{15} & \textbf{20} & \textbf{25} & \textbf{30} & \textbf{35} & \textbf{40} & \textbf{45} & \textbf{50} & \textbf{55} & \textbf{60} \\
\hline
RAGC & 13587.25 & 20370.42 &  27910.27 &  35267.94 &  41522.85 & 48274.80 & 55476.61 & 61705.92 & 66863.54 & 73672.22 & 77190.52\\
RATC & 13088.72 & 19753.98 &  26618.28 &  33074.21 &  39414.35 & 45023.32 & 52048.15 & 57201.57 & 62457.77 & 68722.41 & 71858.87\\
RATS & 11831.57 & 17262.82 &  23973.03 &  28041.96 &  32976.93 & 36313.14 & 40361.33 & 44237.78 & 46178.09 & 48862.71 & 50298.95\\
PGAC & 14267.05 & 21962.87 &  28309.40 &  35157.00 &  41658.65 & 47916.49 & 54059.09 & 58254.65 & 64978.22 & 70702.71 & 74167.68\\
EPSR & 13473.22 &19921.45   & 25287.26 &  31715.40 &  37717.50 & 42347.48 & 48413.21 & 52016.66 & 58032.95 & 62977.44 & 66111.87\\
IDS  & 13718.35 & 20854.14 &  27022.15 &  33184.55 &  39588.14 & 43586.59 & 49353.74 & 54019.48 & 57723.58 & 62196.40 & 64446.65\\
IDA  & 13718.35 & 20854.14 &  27022.15 &  33184.55 &  39588.14 & 43586.59 & 49353.74 & 54019.48 & 57723.58 & 62196.40 & 64446.65\\
\hline
\end{tabular}}
\end{table}
Table~\ref{tab:effective_area} and Figure~\ref{fig:effective_area} compare the effective sensing area obtained by different orientation strategies as the number of deployed sensors increases. The effective sensing area increases monotonically for all methods because additional sensors enlarge the observable region. However, the proposed RAGC consistently achieves the largest effective sensing area over the entire deployment range by explicitly maximizing the visible sensing region while accounting for obstacle-induced occlusions. RATC and RATS achieve slightly smaller sensing areas because their optimization objective progressively shifts from maximizing overall area to improving target visibility and reducing redundant sensor activation through sleep scheduling. Nevertheless, the proposed methods consistently outperform the existing approaches, demonstrating that obstacle-aware robust orientation optimization significantly improves the actual observable sensing region.

\begin{figure}
\begin{subfigure}{0.5\textwidth}
  \centering
  \includegraphics[width=\linewidth]{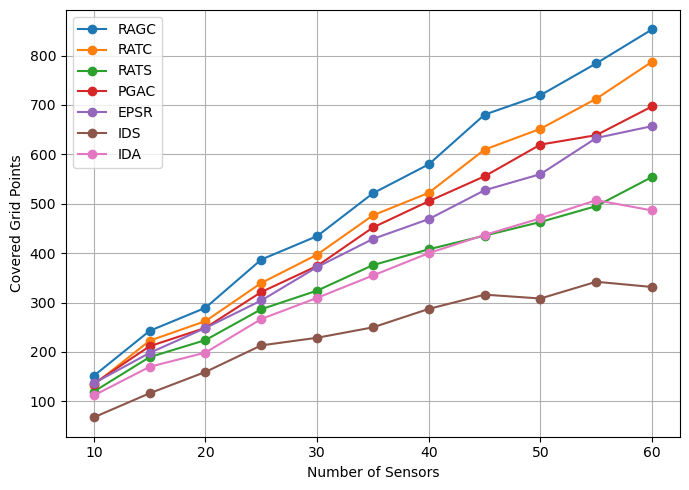}
  \caption{Grid coverage.}
  \label{fig:covered_vs_failed}
\end{subfigure}
\begin{subfigure}{0.5\textwidth}
  \centering
  \includegraphics[width=\linewidth]{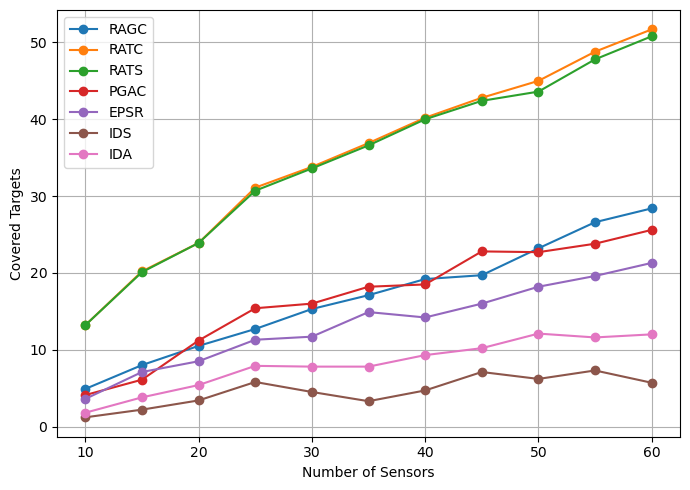}
  \caption{Target coverage.}
  \label{fig:targets_vs_sensors}
\end{subfigure}
\caption{Performance comparison of different optimization methods under varying numbers of deployed sensors in terms of (a) effective grid coverage and (b) target coverage.}
\label{fig:coverage_sensor_variation}
\end{figure}

Figure~\ref{fig:coverage_sensor_variation} illustrates the variation of effective grid coverage and target coverage with respect to the number of deployed sensors. As expected, increasing the sensor density consistently improves both coverage measures for all approaches. However, the proposed methods exhibit a clear advantage over the existing techniques. RAGC achieves the highest grid coverage throughout the entire range of sensor deployments by explicitly maximizing obstacle-aware robust sensing coverage. In contrast, RATC and RATS prioritize target monitoring, resulting in significantly higher target coverage while maintaining competitive grid coverage. Although RATS activates fewer sensors through the proposed sleep scheduling mechanism, its target coverage remains almost identical to RATC, demonstrating that substantial energy savings can be achieved with only negligible performance degradation.

A similar trend is observed when varying other simulation parameters, including the sensing field-of-view, the number of target points and the number of obstacles. The quantitative results are summarized in Table~\ref{tab:parameter_analysis_for_grid_points} and Table~\ref{tab:parameter_analysis_for_targets}, respectively for covered points and covered targets. Increasing the field-of-view generally enlarges the observable region and improves overall coverage, whereas increasing the number of targets provides more opportunities for successful target monitoring. Conversely, a larger number of obstacles reduces the effective sensing region because of visibility occlusions, thereby decreasing achievable coverage for all methods. Nevertheless, the proposed robust optimization framework effectively mitigates these degradations by adapting sensor orientations according to the available obstacle-free regions.
\begin{table}
\centering
\caption{Covered grid points comparison of different optimization strategies under varying network parameters.}
\begin{tabular}{llllllllllll}
\hline
\text{Angle of view} & RAGC & RATC & RATS &  PGAC & EPSR & IDS & IDA\\
\hline 
30 & 404.5 & 358.0 & 255.8 & 309.2 & 310.3 & 192.6 & 247.1 \\
45 & 578.7 & 518.3 & 367.4 & 491.1 & 448.9 & 258.5 & 369.9  \\
60 & 745.3 & 671.8 & 483.5 & 624.8 & 594.1 & 316.1 & 497.6 \\
90 & 1001.5 & 935.8 & 640.6 & 780.8 & 799.1 & 406.7 &583.7 \\
160 & 1429.5 & 1387.5 & 987.4  & 1197.0 & 1181.5 & 572.1  & 928.1   \\
250 & 1641.5 & 1636.1 & 1225.1 & 1385.8 & 1500.2 & 623.9  & 935.3   \\
360& 1657.8 & 1657.5 & 1414.1 & 1543.6 & 1716.0 &655.6   & 1078.8  \\
\hline
\text{Number of targets} & RAGC & RATC & RATS &  PGAC & EPSR & IDS & IDA\\
\hline  
50 & 1621.6 & 1621.3 & 1082.6 & 1385.0 & 1667.2 & 591.6  & 970.7  \\
60 & 1656.2 & 1655.7 & 1248.1 & 1459.1 & 1718.0 & 682.2  & 983.6   \\
70 & 1644.4 & 1644.1  & 1291.6  & 1418.9  & 1699.2 & 655.3  & 951.2 \\
80 & 1622.7 & 1620.2 & 1326.1 & 1439.2 & 1661.2 & 649.3  & 1046.5\\
90 & 1695.7 & 1695.2 & 1470.8 & 1559.4 & 1751.8 & 674.1  & 1192.0 \\
100 & 1695.2 & 1694.7 & 1447.3 & 1519.8 & 1740.3 & 661.0  & 1119.1  \\
150 & 1616.7 & 1616.7 & 1530.8 & 1572.9 & 1657.8 &  633.4 & 1059.1  \\
\hline
\text{Number of obstacles} & RAGC & RATC & RATS &  PGAC & EPSR & IDS & IDA\\
\hline  
0 & 1882.7 & 1882.7 & 1734.6 & 1803.9 & 1882.7 & 720.3  & 1386.0 \\
2 & 1827.8 & 1827.8 & 1708.8 & 1773.9 & 1845.2 & 733.4  & 1346.9  \\
5 & 1631.1 & 1631.1 & 1553.1 & 1583.0 & 1669.7 & 644.2  & 1081.1 \\
7 & 1589.1 & 1584.7 & 1455.3 & 1527.5 & 1656.6 & 660.2  & 1080.3\\
10 & 1468.1 & 1463.2 & 1375.8 & 1425.7 & 1576.3 & 614.9  & 1153.2  \\
\hline
\end{tabular}
\label{tab:parameter_analysis_for_grid_points}
\end{table}

\begin{table}
\centering
\caption{Covered targets comparison of different optimization strategies under varying network parameters.}
\begin{tabular}{llllllllllll}
\hline
\text{Angle of view} & RAGC & RATC & RATS &  PGAC & EPSR & IDS & IDA\\
\hline
30 & 15.0 & 40.6 & 40.1 & 14.0 & 9.0  & 4.6  & 9.7  \\
45 & 19.0 & 43.5 & 43.1 & 17.1 & 13.7 & 2.7  & 5.0   \\
60 & 23.8 & 47.9 & 46.9 & 25.2 & 18.6 & 1.8  & 13.5 \\
90 & 33.7 & 49.3 & 48.0 & 31.1 & 26.6 & 7.7  & 19.7\\
160 & 47.3 & 55.0 & 50.6 & 46.9 &37.8  & 7.2 & 16.9  \\
250 & 52.8 & 54.6 & 49.9 & 54.5 & 49.3 & 20.2 & 31.3  \\
360& 54.7 & 54.7 &52.9  & 56.2 & 56.2 &11.1  &34.1   \\
\hline
\text{Number of targets} & RAGC & RATC & RATS &  PGAC & EPSR & IDS & IDA\\
\hline
50 & 28.5& 33.8 &   33.8 & 35.0 & 35.0 & 13.4 & 23.3 \\
60 & 39.0& 41.9 &   41.9 & 44.3 & 44.3 &  7.3 & 17.8  \\
70 & 43.8& 47.3 &   47.3 & 49.2 & 49.2 &  9.9 & 21.3 \\
80 & 50.1& 53.0 &   53.0 & 54.9 & 54.9 & 19.6 & 32.3\\
90 & 60.0& 64.6 &   64.6 & 66.4 & 66.4 & 15.6 & 38.3  \\
100 & 65.4& 70.2 &   70.2 & 72.4 & 72.4 & 16.3 & 30.9  \\
150& 99.3 & 101.6 & 101.6   & 105.0 & 105.0 &  41.2 &  67.6  \\
\hline
\text{Number of obstacles} & RAGC & RATC & RATS &  PGAC & EPSR & IDS & IDA\\
\hline
0 & 105.5  & 109.6 &  109.6 & 109.6 & 109.6 &  40.4 & 64.4  \\
2 & 104.0  & 107.3 &  107.3 & 108.4 & 108.4 &  27.7 &   63.2 \\
5 & 101.7  & 103.7 &  103.7 & 106.4 & 106.4 & 30.1  & 61.3  \\
7 & 100.3  & 105.1 &  105.1 & 110.3 & 110.3 & 29.7  & 50.7 \\
10 &  98.0  & 100.9 &  100.9 & 108.5 & 108.5 & 12.0  &  42.2  \\
\hline
\end{tabular}
\label{tab:parameter_analysis_for_targets}
\end{table}

Figure~\ref{fig:covered_vs_failed} illustrates the variation of effective grid coverage under progressive sensor failures. At each time interval, approximately 10\% of the remaining active sensors are randomly deactivated to emulate realistic hardware failures or battery depletion. As expected, the grid coverage gradually decreases for all methods as more sensors fail. Nevertheless, RAGC consistently maintains the highest grid coverage throughout the entire failure process because its orientation optimization explicitly maximizes the overall observable sensing region while simultaneously reducing redundant overlaps between neighbouring sensors. RATC closely follows RAGC since it retains many of the same coverage characteristics while additionally emphasizing target monitoring. Although RATS employs fewer active sensors due to sleep scheduling, its coverage degradation remains gradual, demonstrating that the proposed scheduling strategy removes largely redundant sensors while preserving the most informative sensing directions. In contrast, IDS exhibits the lowest grid coverage because each sensor optimizes only its own Voronoi region without considering the global sensing objective. PGAC, EPSR and IDA provide intermediate performance but experience faster degradation as sensor failures increase.
\begin{figure}[!htbp]
\centering

\begin{subfigure}{0.49\textwidth}
    \centering
    \includegraphics[width=\linewidth]{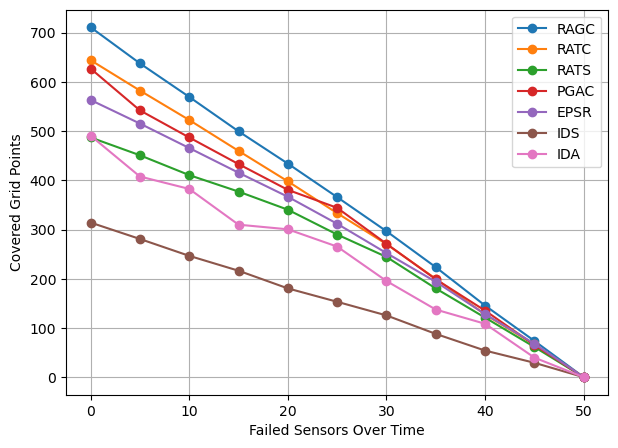}
    \caption{Grid coverage.}
    \label{fig:covered_vs_failed}
\end{subfigure}\hfill
\begin{subfigure}{0.49\textwidth}
    \centering
    \includegraphics[width=\linewidth]{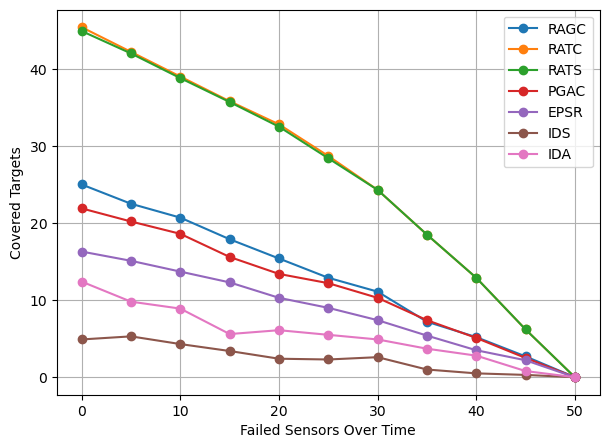}
    \caption{Target coverage.}
    \label{fig:target_vs_failed}
\end{subfigure}

\caption{Comparison of grid and target coverage under progressive sensor failures.}
\label{fig:deployment_optimization}
\end{figure}
Figure~\ref{fig:target_vs_failed} presents the target coverage obtained under the same progressive sensor failure scenario. Unlike grid coverage, the proposed RATC and RATS strategies clearly outperform all competing methods throughout almost the entire experiment. RATC achieves the highest target coverage by explicitly optimizing sensor orientations toward target visibility, whereas RATS maintains nearly identical target coverage despite activating substantially fewer sensors. This demonstrates that the proposed sleep scheduling mechanism successfully identifies redundant sensors without significantly affecting the monitoring capability of the network. Although RAGC provides excellent grid coverage, it covers fewer targets because its optimization objective is distributed over the entire sensing field rather than concentrating on target locations. PGAC and EPSR achieve moderate target coverage, whereas IDS and IDA perform considerably worse due to the absence of explicit target-oriented optimization.

Overall, the experimental results demonstrate that the three proposed strategies satisfy different operational requirements. RAGC is suitable when maximizing overall environmental coverage is the primary objective, RATC is preferable for applications requiring maximum target surveillance and RATS provides the most energy-efficient solution by maintaining almost the same target coverage as RATC while significantly reducing the number of active sensors. Furthermore, the gradual degradation observed under progressive sensor failures confirms that the proposed robust optimization framework remains resilient under uncertain deployment conditions and hardware failures, making it well-suited for long-term aerial directional sensor network applications.

\section{Conclusion and future work}
This paper presents a unified robust optimization framework for aerial directional sensor networks operating in obstacle-rich environments. By integrating obstacle-induced geometric clipping with a grid-based computational representation, the proposed framework establishes a direct connection between realistic visibility modelling and efficient coverage optimization under deployment uncertainty. The RRF is incorporated to ensure robustness against positional uncertainty, while the proposed RAGC, RATC and RATS strategies progressively improve effective grid coverage, target monitoring and energy efficiency through robust orientation optimization and target-aware sleep scheduling. Extensive simulation results demonstrate that the proposed framework consistently outperforms representative approaches under different deployment scenarios, obstacle configurations and sensor failure conditions. Overall, the proposed methodology provides a mathematically rigorous and practically applicable solution for resilient obstacle-aware aerial sensor network design. Future work will investigate dynamic obstacles, heterogeneous sensor networks, adaptive online reorientation and large-scale three-dimensional deployment scenarios.

\section*{Acknowledgements}
\textit{The first author acknowledges the financial support provided by the Indian Institute of Technology Kharagpur, Kharagpur, West Bengal, India-721302.} 

\end{document}